\documentclass[reqno,12pt,a4paper]{amsart}

\usepackage{mathrsfs}
\usepackage{amssymb}
\usepackage{amsfonts}
\usepackage{latexsym}
\usepackage{amsthm}
\usepackage{graphicx}
\def\lb{\label}

\newcommand{\er}[1]{\textrm{(\ref{#1})}}

\begin{document}


\renewcommand{\theequation}{\arabic{section}.\arabic{equation}}
\theoremstyle{plain}
\newtheorem{theorem}{\bf Theorem}[section]
\newtheorem{lemma}[theorem]{\bf Lemma}
\newtheorem{corollary}[theorem]{\bf Corollary}
\newtheorem{proposition}[theorem]{\bf Proposition}
\newtheorem{definition}[theorem]{\bf Definition}
\newtheorem{remark}[theorem]{\bf Remark}

\def\a{\alpha}  \def\cA{{\mathcal A}}     \def\bA{{\bf A}}  \def\mA{{\mathscr A}}
\def\b{\beta}   \def\cB{{\mathcal B}}     \def\bB{{\bf B}}  \def\mB{{\mathscr B}}
\def\g{\gamma}  \def\cC{{\mathcal C}}     \def\bC{{\bf C}}  \def\mC{{\mathscr C}}
\def\G{\Gamma}  \def\cD{{\mathcal D}}     \def\bD{{\bf D}}  \def\mD{{\mathscr D}}
\def\d{\delta}  \def\cE{{\mathcal E}}     \def\bE{{\bf E}}  \def\mE{{\mathscr E}}
\def\D{\Delta}  \def\cF{{\mathcal F}}     \def\bF{{\bf F}}  \def\mF{{\mathscr F}}
\def\c{\chi}    \def\cG{{\mathcal G}}     \def\bG{{\bf G}}  \def\mG{{\mathscr G}}
\def\z{\zeta}   \def\cH{{\mathcal H}}     \def\bH{{\bf H}}  \def\mH{{\mathscr H}}
\def\e{\eta}    \def\cI{{\mathcal I}}     \def\bI{{\bf I}}  \def\mI{{\mathscr I}}
\def\p{\psi}    \def\cJ{{\mathcal J}}     \def\bJ{{\bf J}}  \def\mJ{{\mathscr J}}
\def\vT{\Theta} \def\cK{{\mathcal K}}     \def\bK{{\bf K}}  \def\mK{{\mathscr K}}
\def\k{\kappa}  \def\cL{{\mathcal L}}     \def\bL{{\bf L}}  \def\mL{{\mathscr L}}
\def\l{\lambda} \def\cM{{\mathcal M}}     \def\bM{{\bf M}}  \def\mM{{\mathscr M}}
\def\L{\Lambda} \def\cN{{\mathcal N}}     \def\bN{{\bf N}}  \def\mN{{\mathscr N}}
\def\m{\mu}     \def\cO{{\mathcal O}}     \def\bO{{\bf O}}  \def\mO{{\mathscr O}}
\def\n{\nu}     \def\cP{{\mathcal P}}     \def\bP{{\bf P}}  \def\mP{{\mathscr P}}
\def\r{\rho}    \def\cQ{{\mathcal Q}}     \def\bQ{{\bf Q}}  \def\mQ{{\mathscr Q}}
\def\s{\sigma}  \def\cR{{\mathcal R}}     \def\bR{{\bf R}}  \def\mR{{\mathscr R}}
\def\S{\Sigma}  \def\cS{{\mathcal S}}     \def\bS{{\bf S}}  \def\mS{{\mathscr S}}
\def\t{\tau}    \def\cT{{\mathcal T}}     \def\bT{{\bf T}}  \def\mT{{\mathscr T}}
\def\f{\phi}    \def\cU{{\mathcal U}}     \def\bU{{\bf U}}  \def\mU{{\mathscr U}}
\def\F{\Phi}    \def\cV{{\mathcal V}}     \def\bV{{\bf V}}  \def\mV{{\mathscr V}}
\def\P{\Psi}    \def\cW{{\mathcal W}}     \def\bW{{\bf W}}  \def\mW{{\mathscr W}}
\def\o{\omega}  \def\cX{{\mathcal X}}     \def\bX{{\bf X}}  \def\mX{{\mathscr X}}
\def\x{\xi}     \def\cY{{\mathcal Y}}     \def\bY{{\bf Y}}  \def\mY{{\mathscr Y}}
\def\X{\Xi}     \def\cZ{{\mathcal Z}}     \def\bZ{{\bf Z}}  \def\mZ{{\mathscr Z}}
\def\O{\Omega}

\def\be{{\bf e}}
\def\bv{{\bf v}} \def\bu{{\bf u}}
\def\mn{\mathrm n}
\def\mm{\mathrm m}

\newcommand{\mc}{\mathscr {c}}

\newcommand{\gA}{\mathfrak{A}}          \newcommand{\ga}{\mathfrak{a}}
\newcommand{\gB}{\mathfrak{B}}          \newcommand{\gb}{\mathfrak{b}}
\newcommand{\gC}{\mathfrak{C}}          \newcommand{\gc}{\mathfrak{c}}
\newcommand{\gD}{\mathfrak{D}}          \newcommand{\gd}{\mathfrak{d}}
\newcommand{\gE}{\mathfrak{E}}
\newcommand{\gF}{\mathfrak{F}}           \newcommand{\gf}{\mathfrak{f}}
\newcommand{\gG}{\mathfrak{G}}           
\newcommand{\gH}{\mathfrak{H}}           \newcommand{\gh}{\mathfrak{h}}
\newcommand{\gI}{\mathfrak{I}}           \newcommand{\gi}{\mathfrak{i}}
\newcommand{\gJ}{\mathfrak{J}}           \newcommand{\gj}{\mathfrak{j}}
\newcommand{\gK}{\mathfrak{K}}            \newcommand{\gk}{\mathfrak{k}}
\newcommand{\gL}{\mathfrak{L}}            \newcommand{\gl}{\mathfrak{l}}
\newcommand{\gM}{\mathfrak{M}}            \newcommand{\gm}{\mathfrak{m}}
\newcommand{\gN}{\mathfrak{N}}            \newcommand{\gn}{\mathfrak{n}}
\newcommand{\gO}{\mathfrak{O}}
\newcommand{\gP}{\mathfrak{P}}             \newcommand{\gp}{\mathfrak{p}}
\newcommand{\gQ}{\mathfrak{Q}}             \newcommand{\gq}{\mathfrak{q}}
\newcommand{\gR}{\mathfrak{R}}             \newcommand{\gr}{\mathfrak{r}}
\newcommand{\gS}{\mathfrak{S}}              \newcommand{\gs}{\mathfrak{s}}
\newcommand{\gT}{\mathfrak{T}}             \newcommand{\gt}{\mathfrak{t}}
\newcommand{\gU}{\mathfrak{U}}             \newcommand{\gu}{\mathfrak{u}}
\newcommand{\gV}{\mathfrak{V}}             \newcommand{\gv}{\mathfrak{v}}
\newcommand{\gW}{\mathfrak{W}}             \newcommand{\gw}{\mathfrak{w}}
\newcommand{\gX}{\mathfrak{X}}               \newcommand{\gx}{\mathfrak{x}}
\newcommand{\gY}{\mathfrak{Y}}              \newcommand{\gy}{\mathfrak{y}}
\newcommand{\gZ}{\mathfrak{Z}}             \newcommand{\gz}{\mathfrak{z}}

\def\ve{\varepsilon}   \def\vt{\vartheta}    \def\vp{\varphi}    \def\vk{\varkappa}

\def\A{{\mathbb A}} \def\B{{\mathbb B}} \def\C{{\mathbb C}}
\def\dD{{\mathbb D}} \def\E{{\mathbb E}} \def\dF{{\mathbb F}} \def\dG{{\mathbb G}} \def\H{{\mathbb H}}\def\I{{\mathbb I}} \def\J{{\mathbb J}} \def\K{{\mathbb K}} \def\dL{{\mathbb L}}\def\M{{\mathbb M}} \def\N{{\mathbb N}} \def\dO{{\mathbb O}} \def\dP{{\mathbb P}} \def\R{{\mathbb R}}\def\S{{\mathbb S}} \def\T{{\mathbb T}} \def\U{{\mathbb U}} \def\V{{\mathbb V}}\def\W{{\mathbb W}} \def\X{{\mathbb X}} \def\Y{{\mathbb Y}} \def\Z{{\mathbb Z}}


\def\la{\leftarrow}              \def\ra{\rightarrow}    \def\Ra{\Rightarrow}
\def\ua{\uparrow}                \def\da{\downarrow}
\def\lra{\leftrightarrow}        \def\Lra{\Leftrightarrow}


\def\lt{\biggl}                  \def\rt{\biggr}
\def\ol{\overline}               \def\wt{\widetilde}
\def\no{\noindent}


\let\ge\geqslant                 \let\le\leqslant
\def\lan{\langle}                \def\ran{\rangle}
\def\/{\over}                    \def\iy{\infty}
\def\sm{\setminus}               \def\es{\emptyset}
\def\ss{\subset}                 \def\ts{\times}
\def\pa{\partial}                \def\os{\oplus}
\def\om{\ominus}                 \def\ev{\equiv}
\def\iint{\int\!\!\!\int}        \def\iintt{\mathop{\int\!\!\int\!\!\dots\!\!\int}\limits}
\def\el2{\ell^{\,2}}             \def\1{1\!\!1}
\def\sh{\sharp}
\def\wh{\widehat}
\def\bs{\backslash}
\def\intl{\int\limits}

\def\na{\mathop{\mathrm{\nabla}}\nolimits}
\def\sh{\mathop{\mathrm{sh}}\nolimits}
\def\ch{\mathop{\mathrm{ch}}\nolimits}
\def\where{\mathop{\mathrm{where}}\nolimits}
\def\all{\mathop{\mathrm{all}}\nolimits}
\def\as{\mathop{\mathrm{as}}\nolimits}
\def\Area{\mathop{\mathrm{Area}}\nolimits}
\def\arg{\mathop{\mathrm{arg}}\nolimits}
\def\const{\mathop{\mathrm{const}}\nolimits}
\def\det{\mathop{\mathrm{det}}\nolimits}
\def\diag{\mathop{\mathrm{diag}}\nolimits}
\def\diam{\mathop{\mathrm{diam}}\nolimits}
\def\dim{\mathop{\mathrm{dim}}\nolimits}
\def\dist{\mathop{\mathrm{dist}}\nolimits}
\def\Im{\mathop{\mathrm{Im}}\nolimits}
\def\Iso{\mathop{\mathrm{Iso}}\nolimits}
\def\Ker{\mathop{\mathrm{Ker}}\nolimits}
\def\Lip{\mathop{\mathrm{Lip}}\nolimits}
\def\rank{\mathop{\mathrm{rank}}\limits}
\def\Ran{\mathop{\mathrm{Ran}}\nolimits}
\def\Re{\mathop{\mathrm{Re}}\nolimits}
\def\Res{\mathop{\mathrm{Res}}\nolimits}
\def\res{\mathop{\mathrm{res}}\limits}
\def\sign{\mathop{\mathrm{sign}}\nolimits}
\def\span{\mathop{\mathrm{span}}\nolimits}
\def\supp{\mathop{\mathrm{supp}}\nolimits}
\def\Tr{\mathop{\mathrm{Tr}}\nolimits}
\def\BBox{\hspace{1mm}\vrule height6pt width5.5pt depth0pt \hspace{6pt}}


\newcommand\nh[2]{\widehat{#1}\vphantom{#1}^{(#2)}}
\def\dia{\diamond}

\def\Oplus{\bigoplus\nolimits}



\def\qqq{\qquad}
\def\qq{\quad}
\let\ge\geqslant
\let\le\leqslant
\let\geq\geqslant
\let\leq\leqslant
\newcommand{\ca}{\begin{cases}}
\newcommand{\ac}{\end{cases}}
\newcommand{\ma}{\begin{pmatrix}}
\newcommand{\am}{\end{pmatrix}}
\renewcommand{\[}{\begin{equation}}
\renewcommand{\]}{\end{equation}}
\def\eq{\begin{equation}}
\def\qe{\end{equation}}
\def\[{\begin{equation}}
\def\bu{\bullet}
\def\bq{\mathbf q}

\title{Effective masses for Laplacians on periodic graphs}

\date{\today}

\author[Evgeny Korotyaev]{Evgeny Korotyaev}
\address{Mathematical Physics Department, Faculty of Physics, Ulianovskaya 2,
St. Petersburg State University, St. Petersburg, 198904,
\ korotyaev@gmail.com,}
\author[Natalia Saburova]{Natalia Saburova}
\address{Department of Mathematical Analysis, Algebra and Geometry, Institute of
Mathematics, Information and Space Technologies, Uritskogo St. 68,
Northern (Arctic) Federal University, Arkhangelsk, 163002,
 \ n.saburova@gmail.com}


\begin{abstract}
We consider Laplacians on periodic both discrete and metric
equilateral graphs. Their spectrum consists of an absolutely
continuous part  (which is a union of non-degenerate spectral bands)
and flat bands, i.e., eigenvalues of infinite multiplicity. We
estimate effective masses associated with the ends of each spectral
band in terms of geometric parameters of graphs. Moreover, in the
case of the beginning of the spectrum we determine two-sided
estimates of the effective mass in terms of geometric parameters of
graphs. The proof is based on the Floquet theory, the factorization
of fiber operators, the perturbation theory and the relation between
effective masses for Laplacians on discrete and metric graphs,
obtained in our paper.
\end{abstract}

\maketitle


\vskip 0.25cm

\section {Introduction and main results}
\setcounter{equation}{0}

\subsection {Introduction}

The effective mass approximation is a standard approach in physics (e.g. in solid state physics).
 By this approach, the dynamics  of a particle in a periodic medium
 is  replaced  by the dynamics of a model particle in a model
(simple) medium. Roughly speaking, in this approach, a complicated Hamiltonian is replaced by the model Hamiltonian $-{\D\/2m}$\,,
where $\D$ is the Laplacian and $m$ is the so-called effective mass.
Note that in physics there are a lot of results (a few million references) associated with  effective masses. We mention  the pioneer paper
of Luttinger-Kohn  \cite{LK55} and new results  \cite{Bea12}, \cite{BTV14}, \cite{CL12}, \cite{Iea11},
\cite{Pea11},  and  see references therein.

Now we give few definitions. Consider the band function $\l(\vt)$ of some Hamiltonian, where $\vt=(\vt_\a)_{\a=1}^d\in \T^d=\R^d/(2\pi\Z)^d$ is the quasimomentum. Let the extremum of the band function $\l(\vt)$ be located at the point $\vt_0\in\T^d$. Then the function $\l(\vt)$ can be
expanded into the Taylor series
\[
\lb{tem}
\begin{aligned}
& \l(\vt)=\l(\vt_0)+\m(\e\,)+O\big(|\e\,|^3\big), \\
&\textrm{ where } \qq
\m(\e\,)={1\/2}\sum_{\a,\,\b=1}^d\,M_{\a\b}\,\e_\a
\e_\b,\qqq \e=\vt-\vt_0=(\e_\a)_{\a=1}^d.
\end{aligned}
\]
We call the quadratic form $\mu(\e\,)$ \emph{the effective form}.
Here the linear terms vanish, since $\l(\vt)$ has an extremum at the
point $\vt_0$. If $\l(\vt)$ has a minimum (maximum) at the point $\vt_0$, then the effective form $\mu(\e\,)$ is non-negative (non-positive). The matrix $M=\{M_{\a\b}\}_{\a,\b=1}^d$ is given by
\[
M_{\a\b}={\pa^2\l(\vt_0)\/\pa\vt_\a\pa\vt_\b}\,,
\]
and the matrix $m=M^{-1}$ represents a tensor, which is called
\emph{the effective mass tensor}. The components of this tensor
depend on the coordinate system in quasimomentum space.
In fact, the complicated Hamiltonian is  replaced  by
 the model Hamiltonian (the quadratic form)  $\m(\x)={1\/2}(m^{-1}\x,\x)$,  where $\x=-i\na$.

\

There are a lot of papers devoted to effective masses. In the case
of the Schr\"odinger  operator with a periodic potential on the real
line the effective mass tensor is the scalar and is called the
effective mass. We recall results in this case:

1) Each effective mass is not degenerate, see \cite{F78}.

2) For 1D periodic case Firsova \cite{F75} constructed the conformal
mapping, i.e., the analytic continuation of the quasi-momentum. This
conformal mapping in terms of the effective masses was studied in
\cite{KK95}.

3) There are global estimates of effective masses in terms of norms
of potentials, or gap lengths, band lengths in \cite{K97}.

4) There is a solution of the inverse problem in terms of the effective
masses in \cite{KK97}, i.e., to recover the periodic potential by
the sequence of effective masses.

Thus, the effective masses in 1D case  are well understood.

The effective mass tensor for the multi-dimensional Schr\"odinger  operator with a periodic potential is considered by many authors, but still there are a lot of problems in this subject.
The so-called effective mass approximation is important
to study different problems. We mention only few ones:

1)  Discrete spectrum in the gap of the Schr\"odinger operator with
a periodic potential perturbed by a decaying potential. There are a
lot of papers about, see \cite{B96}, \cite{BS91}, \cite{S98} and
references therein.

2) Homogenization theory has been used to study periodic elliptic
operators near spectral band edges, see \cite{A08}, \cite{AP05}, \cite{BN11}, \cite{B04}, \cite{BS04},  \cite{HW11}   and references therein.

3)  Nonlinear waves in periodic media, see
  \cite{Sp06}, \cite{IM10} and references therein.

4) Linear and nonlinear waves in honeycomb media,
 see  \cite{AY12}, \cite{AY13}, \cite{FW14}, \cite{K08}  and references therein.

Firstly, we remark  that   the properties of band functions in the case of the periodic magnetic Schr\"odinger operator are more complicated than without  magnetic  fields,
 and the effective mass tensor can be degenerate at some magnetic field, see \cite{S04}, \cite{S06}.
Secondly, in the case of effective masses for graphs we know only one paper \cite{K08}, where estimates of effective masses for zigzag
nanotubes in magnetic fields were considered.

\medskip

In our paper we consider effective masses for Laplacians on periodic discrete and metric equilateral graphs. We describe now our main goals:

1) to estimate effective masses associated with the ends of each spectral
band in terms of geometric parameters of graphs;

2) at the beginning of the spectrum to obtain a precise expression for the effective masses and a two-sided estimate of the effective mass tensor (as a matrix) in terms of geometric parameters of graphs.

The proof for the discrete case is based on the Floquet theory, the factorization of fiber operators and the perturbation theory. The proof for the metric case  is essentially based on the precise relations between the eigenvalues of fiber metric and discrete
Laplacians determined in \cite{KS15c}. These relations  allow to obtain a simple identities between effective
masses for the discrete and metric cases.

The results about effective masses are important for the spectral theory of Schr\"odinger operators on periodic metric graphs. For example, using the precise relations between the eigenfunctions of fiber metric and discrete Laplacians and the estimate of the effective masses associated with the beginning of the spectrum in \cite{KS15c} we prove that

$\bu $ all eigenfunctions of the Laplacians on a periodic metric graph are uniformly bounded;

$\bu $ for Schr\"odinger operators with real integrable potentials
on periodic metric graphs the wave operators exist and
are complete.  Furthermore, the standard Fredholm
determinant is well-defined  without any modification for any
dimension. Moreover, the determinant is analytic in the upper
half-plane and the corresponding S-matrix satisfies the Birman-Krein
identity;

$\bu $ the difference of the resolvents for a Schr\"odinger operator with a real integrable potential and for the corresponding Laplace operator on a periodic metric graph belongs to the trace class for any dimension. Note that for Schr\"odinger operators on $\R^d$, $d\geq2$, this is not holds true.

\subsection{Periodic graphs} There are a lot of results about Laplacians
and Schr\"odinger operators on periodic discrete and metric graphs, see
\cite{BK13}, \cite{HN09}, \cite{HS04}, \cite{KS14} -- \cite{KS15c}, \cite{LP08}, \cite{P12}, \cite{RR07}, \cite{S90}, \cite{SS92} and
references therein. For the case of periodic
graphs we know few papers about estimates of bands and gaps:

(1) Lled\'o and Post \cite{LP08} estimated the positions of spectral bands of  Laplacians both on metric and discrete graphs in terms of eigenvalues of the operator on finite graphs (the so-called eigenvalue bracketing).

(2) Korotyaev and Saburova \cite{KS15a}
described a localization of spectral bands and estimated the Lebesgue measure of the spectrum of Schr\"odinger operators with periodic potentials on periodic discrete graphs in terms of eigenvalues of Dirichlet and Neumann operators on a fundamental domain of the periodic graph.

(3) Korotyaev and Saburova \cite{KS14} considered Schr\"odinger
operators with periodic potentials on periodic discrete graphs and estimated the Lebesgue measure of their spectrum  in terms of geometric parameters of the
graph only. Moreover, they estimated a global variation of the Lebesgue measure of the spectrum and a global variation of gap-length in terms of potentials and geometric parameters of the graph.

(4) Korotyaev and Saburova \cite{KS15b} considered Laplacians on periodic equilateral metric graphs and estimated the Lebesgue measure of the bands and gaps on a finite interval in terms of geometric parameters of
the graph.


\medskip

Let $\G=(V,\cE)$ be a connected infinite graph, possibly  having
loops and multiple edges, where $V$ is the set of its vertices and
$\cE$ is the set of its unoriented edges. An edge connecting vertices
$u$ and $v$ from $V$ will be denoted as the unordered pair
$(u,v)_e\in\cE$ and is said to be \emph{incident} to the vertices.
Vertices $u,v\in V$ will be called \emph{adjacent} and denoted by
$u\sim v$, if $(u,v)_e\in \cE$. We define the degree ${\vk}_v$ of
the vertex $v\in V$ as the number of all its incident edges from
$\cE$ (here a loop is counted twice).

Below we consider locally finite
$\Z^d$-periodic graphs $\G$, $d\geq2$, i.e., graphs satisfying the
following conditions:

{\it 1) $\G$ is equipped with an action of the free abelian group $\Z^d$;

2) the degree of each vertex is finite;

3) the quotient graph  $\G_*=\G/{\Z}^d$ is finite.}

\medskip

\no \textbf{Remark.} 1) We do not assume the graph be embedded into
an Euclidean space. But in main applications such a natural
embedding exists (e.g., in modeling waves in thin branching
"graph-like"\, structures: narrow waveguides, quantum wires,
photonic crystal, blood vessels, lungs, see \cite{BK13},
\cite{P12}). The tight-binding approximation is commonly used to
describe the electronic properties of real crystalline structures
(see, e.g., \cite{A76}). This is equivalent to modeling the material
as a discrete graph consisting of vertices (points representing
positions of atoms) and edges (representing chemical bonding of
atoms), by ignoring the physical characters of atoms and bonds that
may be different from one another, see \cite{S13}. The model gives
good qualitative results in many cases. In this case a simple
geometric model is a graph $\G$ embedded into $\R^d$ in such a way
that it is invariant with respect to the shifts by integer vectors
$m\in\Z^d$, which produce an action of $\Z^d$.

2) We also call the quotient graph
$\G_*=\G/{\Z}^d$ \emph{the fundamental graph} of the periodic graph $\G$.
If $\G$ is embedded into the space $\R^d$
the fundamental graph $\G_*$ is a graph on the surface $\R^d/\Z^d$. The fundamental graph $\G_*=(V_*,\cE_*)$ has the vertex set $V_*=V/\Z^d$ and
the set $\cE_*=\cE/\Z^d$ of unoriented edges.

\medskip

Let $\ell^2(V)$ be the
Hilbert space of all square summable functions $f:V\to \C$, equipped
with the norm
$$
\|f\|^2_{\ell^2(V)}=\sum_{v\in V}|f(v)|^2<\infty.
$$
We define the discrete normalized Laplacian (i.e., the Laplace operator) $\D$ on $f\in\ell^2(V)$ by
\[
\lb{DOLN}
 \big(\D f\big)(v)=
 \frac1{\sqrt{\vk_v}}\sum\limits_{(v,\,u)_e\in\cE}
 \bigg({f(v)\/\sqrt{\vk_v}}-{f(u)\/\sqrt{\vk_u}}\bigg),
 \qquad v\in V,
\]
where ${\vk}_v$ is the degree of the vertex $v\in V$ and all loops in the sum
\er{DOLN} are counted twice.

\subsection{The definition of edge indices.} In order to define the Floquet-Bloch
decomposition \er{raz} of discrete Laplacians we need to introduce
two oriented edges $(u,v)$ and $(v,u)$ for each
unoriented edge $(u,v)_e\in \cE$: the oriented edge starting at
$u\in V$ and ending at $v\in V$ will be denoted as the ordered pair
$(u,v)$. Let $\cA$ and $\cA_\ast$ be the sets of all oriented edges of the periodic graph $\G$ and its fundamental graph $\G_*$, respectively.

We define two surjections
\[
\gf_V:V\rightarrow V_*=V/\Z^d, \qqq \gf_\cA:\cA\rightarrow\cA_*=\cA/\Z^d,
\]
which map each element to its equivalence class.

We use {\it an edge index}, which was introduced in \cite{KS14}. The indices are important to study the spectrum of Laplacians and Schr\"odinger operators on periodic graphs, since the fiber operator is expressed in terms of edge indices of the fundamental graph (see \er{l2.15N}). The estimates of the effective masses are also obtained in terms of edge indices.

Let $\nu=\#V_*$, where $\#A$ is the number of
elements of the set $A$.
We fix any $\nu$ vertices of the periodic graph $\G$, which are not $\Z^d$-equivalent to each other and denote this vertex set by $V_0$.
We will call $V_0$ \emph{a fundamental vertex set of $\G$.}
For any
$v\in V$ the following unique representation holds true:
\[
\lb{Dv} v=v_0+[v], \qquad v_0\in
V_0,\qquad [v]\in\Z^d,
\]
where $v_0+[v]$ denotes the action of $[v]\in\Z^d$ on $v_0\in
V_0$.
In other words, each vertex $v$ can be obtained from a vertex $v_0\in V_0$ by the shift by a vector $[v]\in \Z^d$. We will call $[v]$ \emph{coordinates of the vertex $v$ with respect to the fundamental vertex set $V_0$}.
For any
oriented edge $\be=(u,v)\in\cA$ we define {\bf the edge "index"}
$\t({\bf e})$ as the integer vector by
\[
\lb{in}
\t({\bf e})=[v]-[u]\in\Z^d,
\]
where, due to \er{Dv}, we have
$$
u=u_0+[u],\qquad v=v_0+[v], \qquad u_0,v_0\in V_0,\qquad [u],[v]\in\Z^d.
$$
In general, edge indices  depend on the choice of the set $V_0$.
We note that edges connecting vertices from the fundamental vertex set $V_0$ have zero indices. Edges with nonzero indices will be called {\bf bridges}. The bridges provide the connectivity of the periodic graph.

If $\be$ is an oriented edge of the graph $\G$, then by the
definition of the fundamental graph there is an oriented edge
$\be_*=\gf_{\cA}(\be)$ on $\G_\ast$. For the edge
$\be_*\in\cA_\ast$ we define the edge index $\t(\bf e_*)$ by
\[
\lb{inf}
\t(\bf e_*)=\t(\be).
\]
In other words, edge indices of the fundamental graph $\G_\ast$  are
induced by edge indices of the periodic graph $\G$.
An index of a
fundamental graph edge with respect to the fixed fundamental vertex set $V_0$ is uniquely determined by \er{inf}, since
$$
\t(\be+m)=\t(\be),\qqq \forall\, (\be,m)\in\cA \ts \Z^d.
$$

\subsection{Spectrum of discrete Laplacians}
The discrete Laplacian $\D$ on $\ell^2(V)$ is self-adjoint and has the standard decomposition into a
constant fiber direct integral by
\[
\lb{raz}
\ell^2(V)={1\/(2\pi)^d}\int^\oplus_{\T^d}\ell^2(V_\ast)\,d\vt ,\qqq
U\D U^{-1}={1\/(2\pi)^d}\int^\oplus_{\T^d}\D(\vt)d\vt,
\]
$\T^d=\R^d/(2\pi\Z)^d$,
for some unitary operator $U$. Here $\ell^2(V_\ast)=\C^\nu$ is the fiber space, $\nu=\#V_*$. The precise expression of the
Floquet $\nu\ts\nu$ matrix $\D(\vt)$ for the
Laplacian $\D$ is given by \er{l2.15N}. Note that $\D(0)$ is the Laplacian on $\G_*$.

It is convenient to separate all flat bands  from other bands of the Laplacian.
Recall that  $\l_*$ is an eigenvalue
of $\D$ iff $\l_*$ is an eigenvalue of $\D(\vt)$ for any $\vt\in\T^d$
(see Proposition 4.2 in \cite{HN09}).
Thus, if the operator $\D$ has $r\ge 0$ flat bands, then we denote them by
\[
\lb{fb1}
\l_{j}(\vt)=\const, \qqq \forall\, \n-r< j\le \n,\qq \forall\, \vt\in \T^d.
\]
All other eigenvalues (band functions)
$\l_1(\vt),\dots,\l_{\n-r}(\vt)$ are not constant. They can be
enumerated in non-decreasing order (counting multiplicities) by
\[
\label{eq.3.1} \l_1(\vt )\leq\l_2(\vt
)\leq\ldots\leq\l_{\nu-r}(\vt), \qqq \forall\,\vt\in\T^d.
\]
Note that each $\l_n(\vt)\in [0,2]$, $n\in\N_\n=\{1,\ldots,\n\}$. Since $\D(\vt)$
is self-adjoint and analytic in $\vt\in\T^d$, each band function
$\l_n(\cdot)$, $n\in\N_\n$, is a real and piecewise analytic
function on the torus $\T^d$ and creates the spectral band
$\s_n(\D)$ given by
\[
\lb{ban.1} \s_n(\D)=[\l_n^-,\l_n^+]=\l_n(\T^d)\ss [0,2].
\]
Thus, the spectrum of the Laplacian
$\D$ on the discrete periodic graph $\G$ has the form
\[
\lb{r0}
\begin{aligned}
\s(\D)=\bigcup_{\vt\in\T^d}\s\big(\D(\vt)\big)=\bigcup_{n=1}^{\nu}\s_n(\D)=
\s_{ac}(\D)\cup \s_{fb}(\D),\\ \s_{ac}(\D)=\bigcup_{n=1}^{\nu-r}\s_n(\D),\qqq \s_{fb}(\D)=\bigcup_{n=\nu-r+1}^{\nu}\s_n(\D).
\end{aligned}
\]
Here $\s_{ac}(\D)$ is the absolutely continuous spectrum, which is a
union  of non-degenerate intervals from \er{ban.1}, and $\s_{fb}(\D)$ is the set of
all flat bands (eigenvalues of infinite multiplicity). An open
interval between two neighboring non-degenerate spectral bands is
called a \emph{spectral gap}.

\subsection{Effective masses for discrete Laplacians}
Let $\l(\vt)$, $\vt\in\T^d$, be a band function of the Laplacian
$\D$ and let $\l(\vt)$ have a minimum (maximum) at some point
$\vt_0$. Assume that $\l(\vt_0)$ is a simple eigenvalue of
$\D(\vt_0)$. Then the eigenvalue $\l(\vt)$ has the
Taylor series as $\vt=\vt_0+\ve\o$, $\o\in\S^{d-1}$, $\ve=|\vt-\vt_0|\to0$:
\[
\label{lam}
\l(\vt)=\l(\vt_0)+\ve^2\m(\o)+O(\ve^3),\qqq
\textstyle\m(\o)={1\/2}\,\ddot\l(\vt_0+\ve\o)\big|_{\ve=0},
\]
where $\dot u=\pa u/\pa \ve$ and $\S^{d}$ is the $d$-dimensional
sphere.

Now we estimate the effective forms $\m(\o)$ associated with the
ends of each spectral band.

\begin{theorem}\lb{T2}
Let a band function $\l(\vt)$, $\vt\in\T^d$, have a minimum
(maximum) at some point $\vt_0$ and let $\l(\vt_0)$ be a simple
eigenvalue of $\D(\vt_0)$. Then the effective form $\mu(\o)$ from
\er{lam} satisfies
\[
\lb{em11} \big|\mu(\o)\big|\leq {T_1^2\/\r}+T_2\qqq \forall \
\o\in\S^{d-1},
\]
\[
\text{where}\qqq T_s={1\/s}\,\max_{u\in
V_\ast}\sum_{\be=(u,v)\in\cA_\ast} {\|\t
(\be)\|^s\/\sqrt{\vk_u\vk_v}}\,,\qq s=1,2,
\]
where $\r=\r(\vt_0)$ is the distance between $\l(\vt_0)$  and the set
$\s\big(\D(\vt_0)\big)\setminus\big\{\l(\vt_0)\big\}$, and $\t
(\be)$ is the index of the edge $\be\in\cA_\ast$.
\end{theorem}

In order to discuss the effective mass  in the case of the beginning
of the spectrum in Theorem~\ref{Tem} we need recall some notions.
We orient the undirected edges of the fundamental
graph $\G_*$ from the set $\cE_*$ in some arbitrary way
and denote by $\cB_*$ the set of all bridges from $\cE_*$. Let $\n_\ast=\#\cE_\ast$ and let $\n_1=\#\cB_*$. We recall that $\n=\#V_*$.
It is well known (see, e.g., \cite{Ch97}) that $\D(0)=\pa^\ast\pa$, where the $\n_\ast\ts\n$ matrix $\pa
=\{\pa_{\be,v}\}_{\be\in\cE_\ast\atop v\in V_\ast}$ is given by
\[
\lb{Nvt0}
\pa_{\be,v}={1\/\sqrt{\vk_v}}
\left\{
\begin{array}{cl}
  1 & \textrm{if $v$ is the terminal vertex of $\be$} \\
  -1  & \textrm{if $v$ is the initial vertex of $\be$} \\
  0 & \textrm{otherwise}
\end{array}\right..
\]
Since a $\Z^d$-periodic graph $\G$ is connected, there exist $d$
bridges $\be_1,\ldots,\be_d\in\cB_*$ with linearly independent
indices $\t(\be_1),\ldots,\t(\be_d)\in\Z^d$. Let $\cT_0$ be the $d\ts d$ nonsingular matrix whose rows are linearly independent indices of the bridges $\be_1,\ldots,\be_d$ and let $\cT_1$ be the $\n_1\ts d$ matrix, $\n_1=\#\cB_*$, whose rows are indices of all bridges $\be_1,\ldots,\be_{\n_1}$ from the set $\cB_*$, i.e.,
\[\lb{cTT}
\cT_0=\left(\begin{array}{c}
\t(\be_1)\\
\ldots \\
\t(\be_d)
\end{array}\right),\qqq \cT_1=\left(\begin{array}{c}
\t(\be_1)\\
\ldots \\
\t(\be_{\n_1})
\end{array}\right).
\]

We denote by $\L_0>0$ the smallest eigenvalue of the matrix $\cT_0^*\cT_0$ and by $\L_1$ the largest eigenvalue of the matrix $\cT_1^*\cT_1$. For each $\t=(\t_1,\ldots,\t_d)\in\Z^d$ we define the vector $\t^+=(|\t_1|,\ldots,|\t_d|)$.

\begin{theorem}
\label{Tem} Let $\mu(\cdot)$ be the effective form at the beginning
of the spectrum defined in \er{lam}. Then each $\mu(\o)$, $\o\in
\S^{d-1}$, satisfies
\[
\lb{em11.1}
\textstyle \mu(\o)=\|P_*h(\o)\|_*^2,
\]
\[
\lb{tseN} 0<{\L_0\/\vk\,\n\,d}\leq\mu(\o)\leq\|h(\o)\|^2_*=
{1\/\vk}\sum\limits_{\be\in\cB_*}\lan\t
(\be),\o\ran^2\leq{\L_1\/\vk}\,,
\]
where $P_*$ is the orthogonal projection of \, $\C^{\n_*}$ onto the kernel of \, $\pa^*$ and the
vector $h(\o)$ is given by
\[
\lb{hhh}
h(\o)=
{i\/\sqrt{\vk}}\big(\lan\t(\be),\o\ran\big)_{\be\in\cE_\ast}\in\C^{\n_\ast},\qqq
\vk=\sum_{v\in V_\ast}\vk_v,
\]
$\n$ is the number of the fundamental graph vertices, $\n_\ast$ is the number of the fundamental graph edges from $\cE_\ast$,
$\|\cdot\|_\ast$ denotes the norm in $\C^{\n_\ast}$, $\t (\be)$ is
the index of the edge $\be\in\cE_\ast$, $\lan\cdot\,,\cdot\ran$ denotes the standard inner product in $\R^d$.

Moreover, the eigenvalues $\L_0$ and $\L_1$ satisfy
\[
\lb{sve}
\left({d-1\/C}\right)^{d-1}\leq\L_0\leq\L_1\leq\max_{\be_0\in\cB_*}\Big\lan\t^+(\be_0),
\sum_{\be\in\cB_*}\t^+(\be) \Big\ran\,,
\]
where $C=\sum\limits_{j=1}^d\|\t (\be_j)\|^2$.
\end{theorem}

\no {\bf Remark.} 1) The upper estimate  for $\mu(\o)$ in \er{tseN} follows directly from \er{em11.1}. The proof of the lower estimate for $\mu(\o)$ is a more complicated problem.
In particular, the lower estimate in \er{tseN} means that the
effective mass tensor associated with the beginning of the spectrum
is finite.

2) For the most popular periodic graphs (the $d$-dimensional square
lattice, the hexagonal lattice, the Kagome lattice, the
face-centered cubic lattice, the body-centered cubic lattice and
etc.) there exist $d$ bridges with indices forming an orthonormal
basis in $\R^d$ (with respect to some fundamental vertex set). In
this case $\L_0=1$ and the lower
estimate of the effective form in \er{tseN} takes the simple form
\[
\lb{stb}
0<{1\/\vk\,\n\,d}\leq\mu(\o).
\]

3) Each Floquet matrix $\D(\vt)$, $\vt\in\T^d$, has the factorization
$\D(\vt)=\na^\ast(\vt)\na(\vt)$, where $\na(\vt)$ is defined in \er{Nvt}.

4) The number of vertices $\n$ and the number of edges $\n_*$ of the
fundamental graph can be arbitrary, but there is the restriction
(see Theorem 1.2 in \cite{KS15a}):
$$
\n+d\leq\n_\ast+1.
$$

5) There are a lot of papers devoted to estimates of the smallest and the largest eigenvalues of a non-negative matrix $A^*A$, where $A$ is an $m\ts n$ real (complex) matrix (see \cite{N07}, \cite{SD97}, \cite{TC07} and references therein).

\

Let a  matrix $M\ge 0$ be defined by the effective form
$\mu(\o)=\lan M\o,\o\ran$. Then  the effective mass tensor has the form
$$
m=M^{-1}.
$$
Let $0\le M_1\le M_2\le\ldots\le M_d$ be the eigenvalues of $M$. Then
the eigenvalues $m_1\ge m_2\ge\ldots\ge m_d$ of $m$ satisfy
$m_s={1\/M_s}$\,, $s=1,\ldots,d$.

\begin{corollary}\label{Tsem}
 Let  $\be_0\in\cB_\ast$.   Then at the beginning  of
the spectrum  the following estimates hold true
\[
\lb{tem1} C_0\leq\Tr M\leq C_1,
\]
\[
\lb{tem2} {\vk\,\n\,d\/\L_0}\ge m_1\ge\ldots\ge m_d\ge {1\/C_1}\,,
\]
where
\[
\lb{tem3} C_0={1\/\vk\,\n}\,\|\t(\be_0)\|^2,\qqq C_1={1\/\vk}
\sum\limits_{\be\in\cB_\ast}\|\t (\be)\|^2.
\]

\end{corollary}

\no {\bf Remark.} 1) There are a lot of physical problems, where
some parameters are expressed in terms of effective masses. For
example, in 3D case there is an effective mass which corresponds to
the density of states effective mass, given by
\[
m_{dens}=\big(\k^2m_1m_2m_3\big)^{1\/3},
\]
where $\k$ is a degeneracy factor (see, e.g., \cite{G90}). Thus, we can
estimate $m_{dens}$ in terms of lattice parameters. In particular,
the estimates \er{tem2} give
\[
{\k^{2\/3}\/C_1}\le m_{dens}\le {\k^{2\/3}\vk\,\n\,d\/\L_0}\,.
\]

2) The number of bridges on $\G_\ast$ and their indices depend on
the choice of the fundamental vertex set $V_0$. In order to get the best estimate in \er{tseN} and \er{tem1} we  have to choose this set $V_0$ such that
the numbers $\L_1$ and $\sum\limits_{\be\in\cB_\ast}\|\t
(\be)\|^2$ are minimal and the number $\L_0$ is maximal.

3) The proof of the lower estimates in \er{tseN} and in \er{tem1} is
based on the inequality
\[
\lb{leef}
{1\/\vk\,\n}\,\lan \t(\be),\o\ran^2 \leq\m(\o),\qqq \forall\,
(\be,\o)\in\cB_\ast\ts \S^{d-1},
\]
which is proved in Proposition \ref{Pcy1}.

\subsection{Examples of effective masses}
We consider the hexagonal lattice $\bG$, shown in Fig.\ref{ff.0.3}\emph{a}. The lattice $\bG$ is invariant
under translations through the vectors $a_1$, $a_2$.
The fundamental vertex set $V_0=\{v_1,v_2\}$ and the fundamental graph $\bG_*$ are also shown in the figure.

\setlength{\unitlength}{1.0mm}
\begin{figure}[h]
\centering

\unitlength 1mm 
\linethickness{0.4pt}
\ifx\plotpoint\undefined\newsavebox{\plotpoint}\fi 
\begin{picture}(100,45)(0,0)

\put(5,10){(\emph{a})}

\put(14,10){\circle{1}}
\put(28,10){\circle{1}}
\put(34,10){\circle{1}}
\put(48,10){\circle{1}}

\put(18,16){\circle{1}}
\put(24,16){\circle*{1}}
\put(38,16){\circle{1}}
\put(44,16){\circle{1}}

\put(14,22){\circle{1}}
\put(28,22){\circle*{1}}
\put(34,22){\circle{1}}
\put(48,22){\circle{1}}

\put(18,28){\circle{1}}
\put(24,28){\circle{1}}
\put(38,28){\circle{1}}
\put(44,28){\circle{1}}

\put(14,34){\circle{1}}
\put(28,34){\circle{1}}
\put(34,34){\circle{1}}
\put(48,34){\circle{1}}

\put(18,40){\circle{1}}
\put(24,40){\circle{1}}
\put(38,40){\circle{1}}
\put(44,40){\circle{1}}

\put(28,10){\line(1,0){6.00}}
\put(18,16){\line(1,0){6.00}}
\put(38,16){\line(1,0){6.00}}

\put(28,22){\line(1,0){6.00}}
\put(18,28){\line(1,0){6.00}}
\put(38,28){\line(1,0){6.00}}

\put(28,34){\line(1,0){6.00}}
\put(18,40){\line(1,0){6.00}}
\put(38,40){\line(1,0){6.00}}

\put(14,10){\line(2,3){4.00}}
\put(34,10){\line(2,3){4.00}}
\put(24,16){\line(2,3){4.00}}
\put(44,16){\line(2,3){4.00}}

\put(14,22){\line(2,3){4.00}}
\put(34,22){\line(2,3){4.00}}
\put(24,28){\line(2,3){4.00}}
\put(44,28){\line(2,3){4.00}}

\put(14,34){\line(2,3){4.00}}
\put(34,34){\line(2,3){4.00}}

\put(28,10){\line(-2,3){4.00}}
\put(48,10){\line(-2,3){4.00}}
\put(38,16){\line(-2,3){4.00}}
\put(18,16){\line(-2,3){4.00}}

\put(28,22){\line(-2,3){4.00}}
\put(48,22){\line(-2,3){4.00}}
\put(38,28){\line(-2,3){4.00}}
\put(18,28){\line(-2,3){4.00}}

\put(28,34){\line(-2,3){4.00}}
\put(48,34){\line(-2,3){4.00}}

\put(30,18){$\scriptstyle a_1$}
\put(20.5,22){$\scriptstyle a_2$}

\put(24,16){\vector(0,1){12.0}}
\put(33,21.3){\vector(3,2){0.5}}

\qbezier(24,16)(29,19)(34,22)

\put(24.8,21.5){$\scriptstyle v_1$}
\put(35,21.5){$\scriptstyle v_2+a_1$}
\put(25.0,28.0){$\scriptstyle v_2+a_2$}
\put(22.0,13.5){$\scriptstyle v_2$}

\put(75,10){\circle*{1}}
\put(83,21){\circle*{1}}

\put(75,30){\circle{1}}
\put(95,40){\circle{1}}
\put(95,20){\circle{1}}

\put(75,10){\vector(0,1){20.0}}
\put(75,10){\vector(2,1){20.0}}

\multiput(95,20)(0,7){3}{\line(0,1){4}}
\put(75,30){\line(2,1){4.0}}
\put(82,33.5){\line(2,1){4.0}}
\put(89,37){\line(2,1){4.0}}

\qbezier(83,21)(89,20.5)(95,20)
\qbezier(83,21)(79,15.5)(75,10)
\qbezier(83,21)(79,25.5)(75,30)

\put(71,8.0){$v_2$}
\put(83,22){$v_1$}
\put(96,19.0){$v_2$}
\put(71.0,31.0){$v_2$}
\put(93.5,42.0){$v_2$}
\put(85,13){$a_1$}
\put(70,20.0){$a_2$}

\put(76,17.2){$\mathbf{e}_1$}
\put(89,21.2){$\mathbf{e}_2$}
\put(78,26.6){$\mathbf{e}_3$}

\put(64,10){(\emph{b})}
\end{picture}

\vspace{-0.5cm} \caption{\footnotesize  \emph{a}) Graphene $\bG$; \quad \emph{b}) the fundamental graph $\bG_\ast$ of the graphene.} \label{ff.0.3}
\end{figure}

The effective form associated with the beginning of the spectrum is given by
\[\lb{hex}
\textstyle\m(\o)={1\/9}\,\big(\o_1^2+\o_2^2-\o_1\o_2\big),\qqq
\forall\, \o=(\o_1,\o_2)\in\S^1
\]
(for more details see subsection 5.2). Since the hexagonal lattice is a bipartite graph, then the effective form associated with the upper end of the spectrum is $-\m(\o)$.
From \er{hex} it follows that
\[
\textstyle{1\/18}\leq\m(\o)\leq{1\/6}\,.
\]
On the other hand, the fundamental graph $\bG_*$ has only two bridges with indices $(1,0)$, $(0,1)$, $d=2$, $\vk=6$, $\n=2$. Then $\L_0=\L_1=1$ and the estimate \er{tseN} gives
\[
\textstyle{1\/24}\leq\mu(\o)\leq{1\/6}\,.
\]
We note that there is no gap in the spectrum of the Laplacian on the hexagonal lattice.

\

\no \textbf{Remark.} 1) In the present paper we consider effective
masses for the both normalized Laplacians $\D$ defined by \er{DOLN} and combinatorial Laplacians $\wh\D$ defined by \er{DOL} (see Section~\ref{Sec4}). In the case of a
regular graph (i.e. the one with constant degrees of vertices), the
spectra of normalized and combinatorial Laplacians can be easily
related. However, these operators are different.  For example, if a
graph is not regular, then these spectra need to be studied
independently. Moreover, for a bipartite graph the spectrum of the
Laplacian $\D$ is symmetric with respect to the point 1. Then
we need to study only a half of the spectral bands. For the
combinatorial Laplacian $\wh\D$ this property holds true only for
bipartite regular graphs (in this case the spectrum is symmetric
with respect to the point $\vk_0$, where $\vk_0$ is the degree of all
vertices of the graph).

2) In Section \ref{Sec3} we also consider effective masses for the
Laplacian on metric graphs.

\section{\lb{Sec2} Proof of the main results}
\setcounter{equation}{0}

\subsection{Floquet decomposition of Laplacians}
In Proposition \ref{pro0} we present a property of edge indices needed to prove Theorem \ref{pro2}.

\begin{proposition}\label{pro0}
Let $\t^{(1)}(\be_*)$ be the index of an edge
$\be_*=(u_*,v_*)\in\cA_*$ with respect to another fundamental vertex set $V_1$ and let $u_1,v_1\in V_1$ such that $u_*=\gf_V(u_1)$, $v_*=\gf_V(v_1)$. Then
\begin{equation}\label{ind'}
\t^{(1)}(\be_*)=\t(\be_*)+[u_1]-[v_1],
\end{equation}
where $[v]$ is the coordinates of the vertex $v$ with respect to the fundamental vertex set $V_0$.
\end{proposition}
\no {\bf Proof.} Let $\be_*=(u_*,v_*)$ be an oriented edge of the fundamental graph with an index $\t(\be_*)$. Then by the definition of the fundamental graph and the formulas \er{in}, \er{inf} there is an edge $\be=(u_0,v_0+\t(\be_*))\in\cA$,
where $u_0,v_0\in V_0$ such that $u_*=\gf_V(u_0)$, $v_*=\gf_V(v_0)$.
Due to \er{Dv}, for the vertex $u_1,v_1\in V_1$ we have
\[\lb{ii1}
u_1=u_0+[u_1],\qqq v_1=v_0+[v_1]
\]
for some $[u_1],[v_1]\in\Z^d$. This yields
\[
u_0=u_1-[u_1],\qq v_0+\t(\be_*)=v_1-[v_1]+\t(\be_*)
\]
and, due to the definition \er{in} of the edge index,
\[
\t^{(1)}(\be_*)=\t^{(1)}(\be)=-[v_1]+\t(\be_*)+[u_1].
\]
\BBox

\begin{theorem}\label{pro2}
i) The Laplacian $\D$ acting on $\ell^2(V)$ has the
decomposition into a constant fiber direct integral \er{raz},
where the Floquet (fiber) matrix $\D(\vt)=\{\D_{uv}(\vt)\}_{u,v\in V_\ast}$ for the Laplacian $\D$ has the form
\[
\lb{l2.15N}
\D_{uv}(\vt)=\d_{uv}-\ca \displaystyle {1\/\sqrt{\vk_u\vk_v}}
\sum\limits_{\be=(u,v)\in\cA_\ast} e^{\,i\lan\t (\be),\,\vt\ran },
\qq &  {\rm if}\  \ (u,v)\in \cA_\ast \\
\qqq 0, &  {\rm if}\  \ (u,v)\notin \cA_\ast \ac,
\]
$\d_{uv}$ is the Kronecker delta, $\vk_v$ is the degree of the vertex $v$, $\t (\be)$ is the index of the edge $\be\in\cA_\ast$,
$\lan\cdot\,,\cdot\ran$ denotes the standard inner
product in $\R^d$.

ii) Let the matrix $\D^{(1)}(\vt)$ be defined by (\ref{l2.15N}) with
respect to another fundamental vertex set $V_1$. Then the matrices
$\D^{(1)}(\vt)$ and $\D(\vt)$ are unitarily equivalent for each $\vt
\in\T^d$.
\end{theorem}
\no {\bf Proof.} i) We omit the proof, since it repeats the proof of Theorem 1.1.i from \cite{KS14}.

ii) The identity \er{ind'} gives that for each $(u,v)\in\cA_\ast$
\begin{equation}
\label{ind1}
\t^{(1)}(u,v)=\t(u,v)+[u_1]-[v_1],
\end{equation}
where $\t^{(1)}(u,v)$ is the index of the edge $(u,v)$ with respect to the fundamental vertex set $V_1$, $u_1,v_1\in V_1$ such that $u=\gf_V(u_1)$, $v=\gf_V(v_1)$, $[v]$ is the coordinates of the vertex $v$ with respect to the fundamental vertex set $V_0$.

Using \er{ind1} we rewrite the entries of $\D^{(1)}(\vt)$ defined by
(\ref{l2.15N}) in the form
\begin{multline*}
\D_{uv}^{(1)}(\vt )=\d_{uv}-\dfrac1{\sqrt{\vk_u\vk_v}}
\sum\limits_{\be=(u,v)\in\cA_*}e^{i\lan\t^{(1)}
(\be),\vt\ran}\\=\d_{uv}-\dfrac{e^{i\lan[u_1]-[v_1],\vt\ran
}}{\sqrt{\vk_u\vk_v}}\sum\limits_{\be=(u,v)\in\cA_*}e^{i\lan\t
(\be),\,\vt\ran } =e^{i\lan[u_1]-[v_1],\,\vt\ran }\D_{uv}(\vt).
\end{multline*}
We define the diagonal $\nu\times\nu$ matrix
$$
\cU(\vt)=\mathrm{diag}\big(
   e^{i\lan [u_1],\,\vt\ran}\big)_{u\in V_*},\qquad \forall\, \vt \in\T^d.
$$
A direct calculation yields
$$
\cU(\vt )\,\D(\vt )\,\cU^{-1}(\vt
)=\D^{(1)}(\vt ),\qqq \forall\, \vt \in\T^d.
$$
Thus, for each $\vt \in\T^d$ the matrices $\D^{(1)}(\vt )$ and $\D(\vt )$ are unitarily equivalent. \qq \BBox

\

\subsection{Floquet matrix factorization}
In order to prove our next result we need to choose the fundamental vertex set $V_0$ in some special way. Let $T=(V_T,\cE_T)$ be a subgraph of the periodic graph $\G$ satisfying the following conditions:

1) \emph{$T$ is a tree, i.e., a connected graph without cycles;}

2) {\it $V_T$ is a fundamental vertex set, i.e., $V_T$ consists of $\n$ vertices of $\G$, which are not $\Z^d$-equivalent to each other. Recall that $\nu$ is the number of vertices of $\G_\ast$.}

\

\no \textbf{Remark.} 1) We need to note that such a graph $T$ always exists, since the
periodic graph is connected, and $T$ is not unique.

2) The graph $T_*=T/{\Z}^d$ is a spanning tree of
the fundamental graph $\G_*$, i.e., $T_*=(V_*,\cE_{T_*})$ is
a subgraph of $\G_*$, which has $\n-1$ edges and contains no cycles.

From now on we assume that the fundamental vertex set $V_0$
coincides with the vertex set $V_T$. Then, by the definition of the
edge index, \textbf{all edges of the spanning tree $T_*$ have zero indices}.

\

We recall a simple fact about spanning trees of a connected graph
(see, e.g., Lemma 5.1 in \cite{B74}).

\begin{lemma}\lb{LIB}
Let $\be\in\cA_\ast$ be an edge of the fundamental graph $\G_\ast$, which is not in the tree $T_*$. Then

i) There exists a unique cycle in $\G_\ast$, containing only $\be$
and edges of $T_*$.

ii) The length of this cycle, i.e., the number of its edges,  is not
more than $\n$ and the sum of all indices of the cycle edges is
$\t(\be)$, where $\t(\be)$ is the index of $\be$.
\end{lemma}

\no \textbf{Remark.} Item ii) follows from i) and the facts that the number of  the edges of $T_*$ is $\n-1$ and all edges of $T_*$ have zero indices.

\

Recall that we orient the undirected edges of the fundamental
graph $\G_*$ from the set $\cE_*$ in some arbitrary way and denote
by $\cB_*$ the set of all bridges from $\cE_*$.

\begin{proposition}\lb{TP1}
i) Each Floquet matrix $\D(\vt)$, $\vt\in\T^d$, has the following factorization
\[
\lb{Dvt}
\D(\vt)=\na^\ast(\vt)\na(\vt),
\]
where the $\n_\ast\ts\n$ matrix $\na(\vt)
=\{\na_{\be,v}(\vt)\}_{\be\in\cE_\ast\atop\,v\in V_\ast}$ is given by
\[
\lb{Nvt}
\na_{\be,v}(\vt)={1\/\sqrt{\vk_v}}\ca e^{i\lan\t
(\be),\,\vt\ran }, \qq &  \textrm{if $v$ is the terminal vertex of $\be$} \\
\qq-1, \qq &  \textrm{if $v$ is the initial vertex of $\be$} \\
e^{i\lan\t
(\be),\,\vt\ran }-1, \qq &  \textrm{if $\be$ is a loop at the vertex $v$} \\
\qqq 0, &  \textrm{otherwise}\ac,
\]
$\n_\ast$ is the number of edges of the fundamental graph
$\G_\ast$ from the set $\cE_\ast$, $\nu$ is the number of vertices of $\G_\ast$, $\lan\cdot\,,\cdot\ran$ denotes the standard inner product in $\R^d$.

ii) For each $\vt\in\T^d\sm \{0\}$, the rank of the matrix
$\na(\vt)$ is equal to $\n$. The rank of the matrix $\na(0)$ is
equal to $\n-1$.

\end{proposition}
{\bf Proof.} i) Define the matrix-valued function
$X(\vt)=\{X_{uv}(\vt)\}_{u,v\in V_\ast}=\na^\ast(\vt)\na(\vt)$. We
show that $\D(\vt)=X(\vt)$. Due to the definition \er{Nvt} of the
matrix $\na(\vt)$,  we have for $u\neq v$:
\[\lb{Xuv}
\begin{aligned}
X_{uv}(\vt)=\sum_{\be\in\cE_\ast}\na^\ast_{u,\be}(\vt)\na_{\be,v}(\vt)=
\sum_{\be=(u,v)\in\cE_\ast}\na^\ast_{u,\be}(\vt)\na_{\be,v}(\vt)+
\sum_{\be=(v,u)\in\cE_\ast}\na^\ast_{u,\be}(\vt)\na_{\be,v}(\vt)\\=
-{1\/\sqrt{\vk_u\vk_v}}\sum_{\be=(u,v)\in\cE_\ast}e^{\,i\lan\t
(\be),\,\vt\ran}-
{1\/\sqrt{\vk_u\vk_v}}\sum_{\be=(v,u)\in\cE_\ast}e^{-i\lan\t
(\be),\,\vt\ran}=-{1\/\sqrt{\vk_u\vk_v}}\sum_{\be=(u,v)\in\cA_\ast}e^{\,i\lan\t
(\be),\,\vt\ran},
\end{aligned}
\]
and for $u=v$
\[
\begin{aligned}\lb{Xuu}
X_{uu}(\vt)=\sum_{\be\in\cE_\ast}\na^\ast_{u,\be}(\vt)\na_{\be,u}(\vt)=
1-{1\/\vk_u}\sum_{\be=(u,u)\in\cE_\ast}2+
{1\/\vk_u}\sum_{\be=(u,u)\in\cE_\ast}\big|e^{\,i\lan\t
(\be),\,\vt\ran}-1\big|^2\\=1-{2\/\vk_u}\sum_{\be=(u,u)\in\cE_\ast}\cos\lan\t
(\be),\vt\ran.
\end{aligned}
\]
Comparing \er{Xuv}, \er{Xuu} with \er{l2.15N}, we conclude that $\D(\vt)=X(\vt)$.

ii) Let $\na_{\centerdot\, v}(\vt)$, $v\in V_\ast$, denote the column of the matrix $\na(\vt)$ corresponding to the vertex $v$. In order to show that the columns of this matrix are linearly independent, we consider their linear combination with coefficients $\a_v$, $v\in V_\ast$:
\[\lb{lic}
\sum_{v\in V_\ast}\a_v\,\na_{\centerdot\,v}(\vt)=0.
\]
From this, using the form \er{Nvt} of the matrix $\na(\vt)$, we obtain
\[
\a_u=\a_ve^{i\lan\t(\be),\,\vt\ran},\qqq \forall\,
\be=(u,v)\in\cE_\ast.
\]
Since $\t(\bar\be)=-\t(\be)$, where $\bar\be=(v,u)$, the last formula can be rewritten in the form
\[\lb{aa'}
\a_u=\a_ve^{i\lan\t(\be),\,\vt\ran},\qqq \forall\,
\be=(u,v)\in\cA_\ast.
\]

Let $\wt\be_1,\ldots,\wt\be_d$ be bridges on the fundamental graph $\G_\ast$
with linearly independent indices
$\t(\wt\be_1),\ldots,\t(\wt\be_d)$. Lemma
\ref{LIB} gives that on $\G_\ast$ there exist cycles $W_s$, $s\in\N_d$, each of which contains only $\wt\be_s$ and edges of $T_*$ and
\[\lb{ner}
\sum_{\be\in W_s}\t(\be)=\t(\wt\be_s)\neq0.
\]
For each vertex $v$ of the cycle $W_s$ the identity \er{aa'} and \er{ner} give
\[\lb{cyc'}
\a_v\big(e^{i\lan\t(\wt\be_s),\,\vt\ran}-1\big)=0.
\]

First, we assume that $\a_v\neq0$ for each $v\in V_\ast$. Then, due to \er{cyc'}, we have
\[
\lan\t(\wt\be_s),\,\vt\ran=0 ,\qqq \forall\, s\in\N_d.
\]
This yields  $\vt=0$, since the vectors
$\t(\tilde{\be}_1),\ldots,\t(\tilde{\be}_d)\in\Z^d$ are linearly
independent. The matrix $\na(0)\vk^{1\/2}$, where $\vk=\diag(\vk_v)_{v\in V_*}$, is just the incidence
matrix of $\G_\ast$. It is known (see, e.g., Proposition 4.3 in \cite{B74})
that the incidence matrix of a connected graph with $\n$ vertices
has rank $\n-1$. From this it follows that the matrix $\na(0)$ also has rank $\n-1$.

Second, let there exist a vertex $v\in V_\ast$ such that $\a_v=0$.
Due to the connectivity of $\G_\ast$ and the identity \er{aa'}, all coefficients $\alpha_v$ in \er{lic} are zeroes. Thus, for all $\vt\in\T^d\sm\{0\}$ the columns of the matrix $\na(\vt)$ are linearly independent and $\rank\na(\vt)=\n$. \qq \BBox

\

\no \textbf{Remark.} 1) The precise expressions for the Floquet matrix
\er{l2.15N} and the factorization \er{Dvt} can not be written without
a notion of the edge index $\t(\be)$ introduced in Subsection 1.3.

2) The matrix $\pa$, defined by \er{Nvt0}, is equal to $\na(0)$.

\

The Taylor expansion of $\na(\vt)$ about the point $\vt_0=0$ is given by
\[
\lb{nas}
\na(\vt)=\na_0+\ve\na_1(\o)+\ve^2\,\na_2(\o)+O(\ve^3) \qq
\textrm{as } \ \vt=\ve\o, \qqq \ve=|\vt\,|\rightarrow0,
\]
where
\[
\lb{anb}
\na_0=\na(0),\qq \na_1(\o)={\dot\na(\ve\o)}\Big|_{\ve=0},\qq
\na_2(\o)={1\/2}\,\ddot\na(\ve\o)\Big|_{\ve=0}, \qq\o\in \S^{d-1}\,.
\]

We also recall the well-known facts about the kernel of the operator $\na^*(0)$ and its relationship with properties of the graph $\G_*$ (for more details see, e.g., Section 4 in \cite{B74}). For each cycle  $W$ on $\G_*$ we define a vector
$\xi_W=\big(\xi_W(\be)\big)_{\be\in\cE_\ast}\in\C^{\n_\ast}$ by
\[\lb{xsi}
\xi_W(\be)=\left\{
\begin{array}{rc}
  1, & \textrm{if }  \be\in W\\
  -1, & \textrm{if }  \bar\be\in W\\
  0, & \textrm{otherwise}
\end{array}\right.,
\]
where $\bar\be=(v,u)$ is the inverse edge for $\be=(u,v)$.

\begin{lemma}\lb{TL2}
i) The dimension of the kernel of $\na_0^*=\na^*(0)$ is $\n_*-\n+1$.

ii) Let $W$ be a cycle on $\G_*$. Then the vector  $\xi_W$, defined
by \er{xsi}, belongs to the kernel of the operator $\na_0^*$, i.e.,
$\na_0^*\xi_W=0$.

iii) As $\be$ runs through the set $\cE_*\setminus\cE_{T_*}$, the $\n_*-\n+1$ elements $\xi_{(T,\be)}\in\C^{\n_*}$, defined by \er{xsi}, form a basis for the kernel of $\na_0^*$.
\end{lemma}

\subsection{Estimates of effective masses}
Let $\l(\vt)$ be a band function for $\D(\vt)$. Assume that
$\l(\vt_0)$ is a simple eigenvalue of $\D(\vt_0)$ for some $\vt_0$
with a normalized eigenfunction
$\p(\vt_0,\cdot)\in\C^\n$.
Then the eigenvalue $\l(\vt)$ and the
corresponding normalized eigenfunction $\p(\vt,\cdot)$ have asymptotics as $\vt=\vt_0+\ve\o$, $\o\in\S^{d-1}$, $\ve\to0$:
\[
\label{lamp}
\begin{aligned}
\l(\vt)=\l(\vt_0)+\ve^2\m(\o)+O(\ve^3), \qqq
\p(\vt,\cdot)=\p_0+\ve\p_1+\ve^2\p_2+O(\ve^3),\\
\textstyle\m(\o)={1\/2}\,\ddot\l(\vt_0+\ve\o)\big|_{\ve=0},\qq
\p_0=\p(\vt_0,\cdot),\qq
\p_1=\p_1(\o,\cdot)=\dot\p(\vt_0+\ve\o,\cdot)\big|_{\ve=0},\\
\textstyle\p_2=\p_2(\o,\cdot)={1\/2}\,\ddot\p(\vt_0+\ve\o,\cdot)\big|_{\ve=0},
\end{aligned}
\]
where $\dot u=\pa u/\pa \ve$ and $\S^{d}$ is the $d$-dimensional
sphere. The Floquet
matrix $\D(\vt)$, $\vt\in\T^d$, defined by \er{l2.15N}, can be
represented in the following form:
\[
\label{geq.1}
\D(\vt)-\l(\vt_0)\1_\n=\D_0+\ve\D_1(\o)+\ve^2\D_2(\o)+O(\ve^3), \\
\]
as $\vt=\vt_0+\ve\o$, $\ve\rightarrow0$, $\o\in
\S^{d-1}$, where
\[\lb{DDD}
\begin{aligned}
\D_0=\D(\vt_0)-\l(\vt_0)\1_\n,\qqq  \D_1(\o)=\dot
\D(\vt_0+\ve\o)\Big|_{\ve=0},\qq \D_2(\o)={1\/2}\,\ddot
\D(\vt_0+\ve\o)\Big|_{\ve=0},
\end{aligned}
\]
$\1_\n$ is the identity $\n\ts\n$ matrix. The equation
$\D(\vt)\p(\vt,\cdot)=\l(\vt)\p(\vt,\cdot)$ after substitution
\er{lamp}, \er{geq.1} takes the form
\[\lb{gaas}
\begin{aligned}
\big(\D_0+\ve\D_1(\o)+\ve^2\D_2(\o)+O(\ve^3)\big)
\big(\p_0+\ve\p_1+\ve^2\p_2+O(\ve^3)\big)\\=
\big(\ve^2\m(\o)+O(\ve^3)\big)
\big(\p_0+\ve\p_1+\ve^2\p_2+O(\ve^3)\big),
\end{aligned}
\]
where $\p_0, \p_1, \p_2$ are defined in \er{lamp}. This asymptotics
gives two identities for any $\o\in \S^{d-1}$:
\[\lb{gve1}
\D_1(\o)\p_0+\D_0\p_1=0,
\]
\[\lb{gve2}
\D_2(\o)\p_0+\D_1(\o)\p_1+ \D_0\p_2=\m(\o)\,\p_0.
\]

\

\no \textbf{Proof of Theorem \ref{T2}.} The Floquet matrix
$\D(\vt)$, $\vt\in\T^d$, defined by \er{l2.15N}, can be represented
in the  form \er{geq.1}. Using \er{l2.15N} we obtain that the
entries of the matrices $\D_s(\o)=\{\D^{(s)}_{uv}(\o)\}_{u,v\in
V_\ast}$, $s=1,2$, defined by \er{DDD}, have the form
\[
\lb{Del1}
\D_{uv}^{(1)}(\o)={-i\/\sqrt{\vk_u\vk_v}}\sum\limits_{\be=(u,\,v)\in{\cA}_\ast}\lan\t
(\be),\o\ran\, e^{\,i\lan\t
(\be),\,\vt_0\ran },
\]
\[
\lb{Del2}
\D_{uv}^{(2)}(\o)={1\/2\sqrt{\vk_u\vk_v}}\sum\limits_{\be=(u,\,v)\in{\cA}_\ast}\lan\t
(\be),\o\ran^2\, e^{\,i\lan\t (\be),\,\vt_0\ran },
\]
for any $\o\in \S^{d-1}$. Since $\l(\vt_0)$ is a simple eigenvalue
of $\D(\vt_0)$, the Floquet matrix $\D(\vt)$, the eigenvalue
$\l(\vt)$ and the corresponding normalized eigenfunction
$\p(\vt,\cdot)$ have asymptotics \er{lamp}, \er{geq.1}, where their
coefficients satisfy the identities \er{gve1}, \er{gve2}.

We recall a simple fact that $\D_1(\o)\p_0$ and $\p_0$ are
orthogonal. Indeed, multiplying both sides of \er{gve1} by $\p_0$
and using that $\D_0\p_0=0$, we have $\lan\D_1\p_0,\p_0\ran=0$,
which yields $\D_1(\o)\p_0\perp\p_0$.

Let $P$ be the orthogonal projection onto the subspace of $\ell^2(V_\ast)$ orthogonal to $\p_0$. From \er{gve1} we obtain
\[\lb{p_1}
\p_1=-(P\D_0)^{-1}P\D_1(\o)\p_0.
\]
Multiplying both sides of \er{gve2} by $\p_0$, substituting \er{p_1}
and using that $\D_0\p_0=0$, we have
\[\lb{sec}
\m(\o)=\lan\D_2(\o)\p_0,\p_0\ran-\lan(P\D_0)^{-1}P\D_1(\o)\p_0,\D_1(\o)\p_0\ran.
\]
This yields
\[\lb{sec1}
\begin{aligned}
\big|\m(\o)\big|\le \big|\lan\D_2(\o)\p_0,\p_0\ran\big|+
\big|\lan(P\D_0)^{-1}P\D_1(\o)\p_0,\D_1(\o)\p_0\ran\big|\\
\leq
\|\D_2(\o)\|+\|(P\D_0)^{-1}P\|\cdot\|\D_1(\o)\|^2\leq
\|\D_2(\o)\|+{1\/\r}\,\|\D_1(\o)\|^2,
\end{aligned}
\]
where $\r=\r(\vt_0)$ is the distance between $\l(\vt_0)$ and $\s\big(\D(\vt_0)\big)\setminus\big\{\l(\vt_0)\big\}$. Due to \er{Del1}, \er{Del2}, we have
\[\lb{NN1}
\|\D_1(\o)\|\leq\max_{u\in V_\ast}\sum_{\be=(u,v)\in\cA_\ast}
{\big|\lan\t (\be),\,\o\ran\big|\/\sqrt{\vk_u\vk_v}}\leq \max_{u\in
V_\ast}\sum_{\be=(u,v)\in\cA_\ast}
{\|\t(\be)\|\/\sqrt{\vk_u\vk_v}}=T_1,
\]
\[\lb{NN2}
\|\D_2(\o)\|\leq\max_{u\in V_\ast}\sum_{\be=(u,v)\in\cA_\ast}
{\lan\t(\be),\,\o\ran^2\/2\sqrt{\vk_u\vk_v}}\leq\max_{u\in
V_\ast}\sum_{\be=(u,v)\in\cA_\ast}{\|\t
(\be)\|^2\/2\sqrt{\vk_u\vk_v}}=T_2.
\]
Substituting \er{NN1}, \er{NN2} into \er{sec1}, we obtain \er{em11}. \qq
\BBox

\

Now we consider an effective mass associated with the beginning of the spectrum.
It is known \cite{SS92} that the lower point of the spectrum of the Laplacian $\D$ is $\l_1(0)=0$. Since the fiber operator $\D(0)$ is the Laplacian on the fundamental graph $\G_\ast$ and $\G_\ast$ is finite and connected,
the following statement holds true (see, e.g., \cite{Ch97}).

\emph{The operator $\D(0)$ has a simple eigenvalue $\l(0)=0$.  The corresponding normalized eigenfunction has the form}
\[\lb{psi0}
\p(0,\cdot)\in\R^\n,\qqq \p(0,v)=\sqrt{\vk_v\/\vk},\qqq v\in V_\ast,
\qqq \vk=\sum_{v\in V_\ast}\vk_v.
\]

\begin{proposition}\lb{Pcy1}
Let $\mu(\o)$ be the effective form at the beginning  of
the spectrum defined in \er{lam}. Then

i) $\mu(\o)$ satisfies the
following identity
\[\lb{mu1'}
\m(\o)=\|\na_1(\o)\p_0\|^2_\ast-\|\na_0\p_1\|_\ast^2=
\|\na_1(\o)\p_0+\na_0\p_1\|_\ast^2,
\]
for any $\o\in \S^{d-1}$, where
\[\lb{n1p0N}
\na_1(\o)\p_0=
{i\/\sqrt{\vk}}\big(\lan\t(\be),\o\ran\big)_{\be\in\cE_\ast}\in\C^{\n_\ast},\qqq
\vk=\sum_{v\in V_\ast}\vk_v,
\]
\[\lb{n0p1N}
\na_0\p_1=
\bigg({\p_1(\o,v)\/\sqrt{\vk_v}}-{\p_1(\o,u)\/\sqrt{\vk_u}}\bigg)
_{\be=(u,v)\in\cE_\ast}\in\C^{\n_\ast}.
\]

ii) the effective form
$\mu(\o)$  satisfies
\[\lb{leab}
\m(\o) \geq{1\/\vk\,\n}\,\lan \t(\be),\o\ran^2,\qqq \forall\,
(\be,\o)\in\cB_*\ts\S^{d-1}.
\]

\end{proposition}
\no {\bf Proof.} i) Since $\l(\vt_0)=0$ at $\vt_0=0$ is a simple eigenvalue
of $\D(0)$, the Floquet matrix $\D(\vt)$, the eigenvalue $\l(\vt)$
and the corresponding normalized eigenfunction $\p(\vt,\cdot)$ have
asymptotics \er{lamp}, \er{geq.1} at $\vt_0=0$, where their coefficients
satisfy the identities \er{gve1}, \er{gve2}. Substituting \er{nas}
into the identity $\D(\vt)=\na^*(\vt)\na(\vt)$ at $\vt=\ve\o$ we
obtain
\begin{multline}
\lb{Dna} \D_0=\na_0^\ast\na_0,  \qq \na_0\p_0=0,\qq
\D_1(\o)=\na_0^\ast\na_1(\o)+\na_1^*(\o)\na_0,\\
\D_2(\o)=\na_0^\ast\na_2(\o)+\na_1^\ast(\o)\na_1(\o)+\na_2^*(\o)\na_0,
\end{multline}
where $\D_0$, $\D_s(\o)$ and $\na_0$, $\na_s(\o)$, $s=1,2$, are defined in \er{DDD} and \er{anb}, respectively.

Multiplying both sides of \er{gve1} by $\p_1$ and both sides of
\er{gve2} by $\p_0$ and using the identities \er{Dna}, we have
\[\lb{ve1'}
\lan\na_1(\o)\p_0,\na_0\p_1\ran+ \|\na_0\p_1\|^2_\ast=0,\qqq
\|\na_1(\o)\p_0\|^2_\ast+ \lan\na_0\p_1,\na_1(\o)\p_0\ran =\m(\o).
\]
Then, combining \er{ve1'}, we obtain
\[\lb{mu2}
\m(\o)=\|\na_1(\o)\p_0\|^2_\ast-\|\na_0\p_1\|_\ast^2.
\]
Relations \er{ve1'} also imply another representation for
$\m(\o)$:
\[\lb{mu1}
\m(\o)=\|\na_1(\o)\p_0+\na_0\p_1\|_\ast^2.
\]
From \er{Nvt} we deduce that the matrix
$\na_1(\o)=\{\na^{(1)}_{\be,v}(\o)\}_{\be\in\cE_\ast\atop\,v\in V_\ast}$ is given by
\[
\label{nab1}
\na^{(1)}_{\be,v}(\o)=\ca \displaystyle i\,{\lan\t
(\be),\,\o\ran\/\sqrt{\vk_v}}\,, &  \textrm{if $v$ is the terminal vertex of $\be$} \\[10pt]
\qqq 0, &  \textrm{otherwise}\ac.
\]
Using \er{nab1}, \er{psi0}, \er{Nvt} and $\na_0=\na(0)$, we obtain \er{n1p0N} and
\er{n0p1N}.

ii) Let $\wt\be\in\cB_\ast$ be a bridge of the fundamental graph $\G_\ast$ with the index $\t(\wt\be\,)$. Then, due to Lemma
\ref{LIB}, on $\G_\ast$ there exists a cycle
$W$
such that $\sum\limits_{\be\in W}\t(\be)=\t(\wt\be\,)$ and the length of $W$ is not more than $\n$.
Let $\xi_W\in\C^{\n_\ast}$ be defined by \er{xsi}.
Then, using \er{mu1} and the Cauchy-Schwarz inequality, we have
\begin{multline}\lb{mmo}
\m(\o)\,\|\xi_W\|_\ast^2=
\|\na_1(\o)\p_0+\na_0\p_1\|_\ast^2\cdot\|\xi_W\|_\ast^2\geq
|\lan\na_1(\o)\p_0+\na_0\p_1,\xi_W\ran_\ast|^2\\=
|\lan\na_1(\o)\p_0,\xi_W\ran_\ast+\lan\na_0\p_1,\xi_W\ran_\ast|^2=
|\lan\na_1(\o)\p_0,\xi_W\ran_\ast+\lan\p_1,\na_0^*\xi_W\ran_\ast|^2,
\end{multline}
where $\lan\cdot\,,\cdot\ran_*$ denotes the standard inner product in $\C^{\n_*}$.
Due to Lemma \ref{TL2}.ii, $\na_0^*\xi_W=0$.
Substituting this identity into \er{mmo} and using that
$\|\xi_W\|_\ast^2\leq\n$ and \er{n1p0N}, we obtain
\begin{multline*}
\vk\,\n\m(\o)\geq \vk\big|\lan\na_1(\o)\p_0,\xi_W\ran_\ast\big|^2=
\Big|\sum_{\be\in\cE_\ast}\xi_W(\be)\lan\t(\be),\o\ran\Big|^2=
\Big\lan\sum_{\be\in\cE_\ast}\xi_W(\be)\t(\be),\o\Big\ran^2\\
=\Big\lan\sum_{\be\in W}\t(\be),\o\Big\ran^2=\lan
\t(\wt\be\,),\o\ran^2,
\end{multline*}
which yields \er{leab}. \qq \BBox

\

\no \textbf{Proof of Theorem \ref{Tem}.} Due to \er{mu1'} we have
\[\lb{mu1''}
\m(\o)=\|\na_1(\o)\p_0+\na_0\p_1\|_\ast^2,
\]
where $\na_1(\o)\p_0$ and $\na_0\p_1$ are defined by \er{n1p0N} and \er{n0p1N}.

We will show that
\[\lb{inc}
\na_1(\o)\p_0+\na_0\p_1\in\ker\na_0^*,\qqq \na_0\p_1\in(\ker\na_0^*)^\perp.
\]
The identities \er{gve1} and \er{Dna} give
$$
0=\D_1(\o)\p_0+\D_0\p_1=\na_0^\ast\na_1(\o)\p_0+\na_1^*(\o)\na_0\p_0+
\na_0^\ast\na_0\p_1=\na_0^\ast\big(\na_1(\o)\p_0+\na_0\p_1\big),
$$
which yields the first formula in \er{inc}. Let $f\in\ker\na_0^*$. Then we obtain
\[
0=\lan \na_0^*f,\p_1\ran=\lan f,\na_0\p_1\ran,
\]
which yields the second formula in \er{inc}.

We recall that $\pa$, defined by \er{Nvt0}, is equal to $\na_0$.
Then from \er{inc} it follows that
\[
\na_1(\o)\p_0+\na_0\p_1=P_*\big(\na_1(\o)\p_0+\na_0\p_1\big)=
P_*\na_1(\o)\p_0+P_*\na_0\p_1=P_*\na_1(\o)\p_0.
\]
Substituting this identity into \er{mu1''}, denoting $h(\o)=\na_1(\o)\p_0$ and using \er{n1p0N}, we obtain \er{em11.1}, \er{hhh}.

Now we prove the lower estimate in \er{tseN}. We recall that on the fundamental graph $\G_\ast$ there exist $d$
bridges  $\be_1,\ldots,\be_d\in\cB_\ast$ with linearly independent
indices
\[\lb{lii}
\t(\be_1)=(\t_{11},\ldots,\t_{1d}),\qq \ldots, \qq \t(\be_d)=(\t_{d1},\ldots,\t_{dd}).
\]
Due to Proposition
\ref{Pcy1}.ii, we have
$$
\m(\o)\geq{1\/\vk\,\n}\,\lan \t(\be_s),\o\ran^2,\qqq \forall\, s\in\N_d.
$$
Summing these estimates over all $s$, we get
\[\lb{esm1}
\m(\o)
\geq{1\/\vk\,\n\,d}\sum_{s=1}^d\lan \t(\be_s),\o\ran^2.
\]
Using \er{lii} we rewrite the quadratic form $\sum\limits_{s=1}^d\lan \t(\be_s),\o\ran^2$ in the form
\[
\begin{aligned}\lb{esm2}
\sum_{s=1}^d\lan \t(\be_s),\o\ran^2=\sum_{s=1}^d\bigg(\sum_{i=1}^d
\t_{si}\,\o_i\bigg)^2=\sum_{s=1}^d\sum_{i,j=1}^d\t_{si}\t_{sj}\,\o_i\,\o_j\\=
\sum_{i,j=1}^d\o_i\,\o_j\sum_{s=1}^d\t_{si}\t_{sj}
=\lan \cT_0^*\cT_0\,\o,\o\ran\geq\L_0,
\end{aligned}
\]
where $\o=(\o_1,\ldots,\o_d)$ and the $d\ts d$ matrix $\cT_0$ is defined in \er{cTT}, $\L_0$ is the smallest eigenvalue of $\cT_0^*\cT_0$.
Since the vectors \er{lii} are linearly independent, from \er{esm1} and \er{esm2} we obtain
\[
\m(\o)
\geq{\L_0\/\vk\,\n\,d}>0\,.
\]
Thus, the lower estimate in \er{tseN} has been proved.

Now we prove the upper estimate in \er{tseN}.
From \er{mu2} we deduce
\[\lb{upe}
\m(\o)\leq\|\na_1(\o)\p_0\|^2_\ast.
\]
Substituting \er{n1p0N} into \er{upe} and using that only bridges have nonzero indices, we obtain
\[
\lb{eem}
\m(\o)\leq{1\/\vk}\sum\limits_{\be\in\cB_\ast}\lan\t
(\be),\o\ran^2.
\]
Similarly \er{esm2} we deduce  that
$\sum\limits_{\be\in\cB_\ast}\lan\t (\be),\o\ran^2\leq\L_1$, where $\L_1$ is the largest eigenvalue of $\cT_1^*\cT_1$ and the matrix $\cT_1$ is defined in \er{cTT}. This and \er{eem} give the upper estimate in \er{tseN}.

Finally, we show \er{sve}. The following simple bounds for $\L_0$ and $\L_1$ were given in \cite{SD97} and \cite{N07}, respectively:
\[\lb{eees}
\begin{aligned}
\L_0>\bigg({d-1\/C}\bigg)^{d-1}(\det \cT_0)^2, \qqq C=\sum_{j,k=1}^d |\t_{jk}|^2,\\
\L_1\leq\max_{1\leq i\leq\n_1}\sum_{j=1}^d|\t_{ij}|\sum_{k=1}^{\n_1}|\t_{kj}|=\max_{\be_0\in\cB_*}\Big\lan\t^+(\be_0),
\sum_{\be\in\cB_*}\t^+(\be) \Big\ran\,,
\end{aligned}
\]
where $\t^+(\be)$ is defined before Theorem \ref{Tem} and
$\t(\be_k)=(\t_{k1},\ldots,\t_{kd})$, $k\in\N_{\n_1}$, is the index of the bridge $\be_k\in\cB_*$.
Since all entries of the matrix $\cT_0$ are integer and $\cT_0$ is nonsingular, $(\det \cT_0)^2\geq1$. Then using \er{eees} we obtain \er{sve}.
\qq \BBox

\

\no \textbf{Proof of Corollary \ref{Tsem}.} We show the upper
estimate in \er{tem1}. Let $M$  be a matrix of the effective form
$\mu(\o)$. And let $M_1,\ldots,M_d$ and $\wt\o_1,\ldots,\wt\o_d$ be
the corresponding eigenvalues and orthonormal eigenvectors of $M$.
Then \er{tseN} implies
\[
\lb{Mse}
0<{\L_0\/\vk\,\n\,d}\leq\m(\wt\o_s)=M_s\leq {1\/\vk}\sum\limits_{\be\in\cB_\ast}\lan\t
(\be),\wt\o_s\ran^2.
\]
Summing the last estimates we have
$$
\Tr M=\sum_{s=1}^d M_s\leq
{1\/\vk}\sum_{s=1}^d\sum\limits_{\be\in\cB_\ast}\lan\t
(\be),\wt\o_s\ran^2={1\/\vk}\sum\limits_{\be\in\cB_\ast}\|\t
(\be)\|^2,
$$
which yields the upper estimate in \er{tem1}. Integrating the
inequality \er{leef} over $\o\in\S^{d-1}$, we obtain the lower
estimate in \er{tem1}.

The estimates \er{tem2} follow from \er{tem1} and \er{Mse}.
\qq \BBox

\section{\lb{Sec3} Effective masses for metric Laplacians}
\setcounter{equation}{0}
\subsection{Metric Laplacians}
In order to define metric Laplacians we identify
 each edge $\be$ of $\G$  with the segment $[0,1]$.
This identification introduces a local coordinate $t\in[0,1]$ along
each edge. Thus, we give an orientation on the edge. An edge starting at $u\in V$ and ending at $v\in V$ will be denoted as the
ordered pair $(u,v)\in\cE$.
Note that the
spectrum of Laplacians on metric graphs does not depend on the
orientation of graph edges. For each function $y$ on $\G$ we define
a function $y_{\be}=y\big|_{\be}$, $\be\in\cE$. We identify each
function  $y_{\be}$ on $\be$ with a function on $[0,1]$ by using the
local coordinate $t\in[0,1]$.
Let $L^2(\G)$ be the Hilbert space of all function $y=(y_\be)_{\be\in\cE}$,
where each $y_\be\in L^2(0,1)$, equipped with the norm
$$
\|y\|^2_{L^2(\G)}=\sum_{\be\in\cE}\|y_\be\|^2_{L^2(0,1)}<\infty.
$$
We define the metric Laplacian $\D_M$ on $y=(y_\be)_{\be\in\cE}\in L^2(\G)$ by
$$
(\D_My)_\be=-y''_\be,\qqq
$$
where $(y''_\be)_{\be\in\cE}\in L^2(\G)$ and
$y$ satisfies the so-called Kirchhoff conditions:
\[
\lb{Dom1}
y \textrm{ is continuous on }\G,\qqq
\sum\limits_{\be\in I(v)}(-1)^{\d(\be,v)}\,y_\be'\big(\d(\be,v)\big)=0, \qq \forall v\in V,
\]
\[\lb{tev}
\d(\be,v)=\left\{
\begin{array}{rl}
1, \qq &\textrm{if $v$ is the terminal vertex of the edge $\be$} \\
 0, \qq &\textrm{if $v$ is the initial vertex of the edge $\be$},
\end{array}\right.
\]
$I(v)$ is the set of all edges of $\G$ incident to the vertex $v\in V$.

The metric Laplacian $\D_M$ on $L^2(\G)$ has the
decomposition into a constant fiber direct integral (see Theorem 1.1 in \cite{KS15c})
\[
\lb{Mraz}
\begin{aligned}
& L^2(\G)={1\/(2\pi)^d}\int^\oplus_{\T^d}L^2(\G_\ast)\,d\vt ,\qqq
\mU\D_M \mU^{-1}={1\/(2\pi)^d}\int^\oplus_{\T^d}\D_M(\vt)d\vt,
\end{aligned}
\]
for some unitary operator $\mU$. The precise expression of the
Floquet operator $\D_M(\vt)$ see in \cite{KS15c}.

Each Floquet operator $\D_M(\vt)$, $\vt\in\T^d$, has the discrete spectrum given by
\[
\lb{egv1}
\begin{aligned}
E_k=(\pi k)^2,\qq k\in\N, \qqq  E_{n,j}(\vt)=z_{n,j}^2(\vt), \qq n\in\N_\n,\qq j=0,1,2,\ldots,\\
z_{n,j}(\vt)=\ca
  z_n(\vt)+\pi j, & j \textrm{ is even} \\
  (\pi-z_n(\vt))+\pi j, & j \textrm{ is odd}
\ac, \qq z_n(\vt)=\arccos(1-\l_n(\vt))\in [0,\pi]
\end{aligned}
\]
(for more details, see \cite{KS15c}).
Thus, the spectrum of the Laplacian
$\D_M$ on the metric equilateral periodic graph $\G$ has the form
\[\lb{sDM1}
\begin{aligned}
\s(\D_M)=\s_{ac}(\D_M)\cup \s_{fb}(\D_M),\qqq
\s_{ac}(\D_M)=\bigcup_{(n,j)\in \N_{\n-r}\ts
\N}\s_{n,j-1}(\D_M),\\
\s_{fb}(\D_M)=\s_{fb}^{(1)}(\D_M)\cup\s_{fb}^{(2)}(\D_M),\qqq
\s_{fb}^{(1)}(\D_M)=\{(\pi k)^2: k\in\N\},\\
\s_{fb}^{(2)}(\D_M)=\bigcup_{n=\nu-r+1}^{\nu}\bigcup_{j=0}^{\infty}\s_{n,j}(\D_M),\qqq
\s_{n,j}(\D_M)=[E_{n,j}^-,E_{n,j}^+],\\
E_{n,j}^\pm=\ca
  \big(z_n^\pm+\pi j\big)^2, & \qq j \textrm{ is even} \\
  \big((\pi-z_n^\mp)+\pi j\big)^2, & \qq j \textrm{ is odd}
\ac,\qqq
z_n^\pm=\arccos(1-\l_n^\pm)\in [0,\pi],
\qq n\in\N_{\n}.
\end{aligned}
\]
Here $\s_{ac}(\D_M)$ is the absolutely continuous spectrum, which is a
union  of non-degenerate intervals, and $\s_{fb}(\D_M)$ is the set of
all flat bands (eigenvalues of infinite multiplicity).

\setlength{\unitlength}{1.0mm}
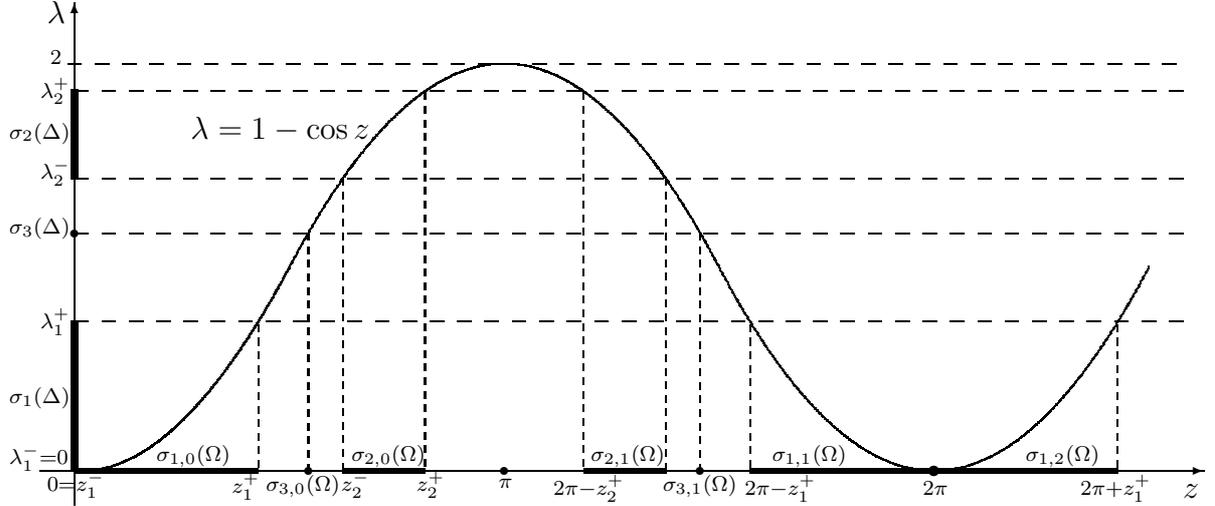
\begin{figure}[h]
\centering
\unitlength 0.9mm 
\linethickness{0.4pt}
\ifx\plotpoint\undefined\newsavebox{\plotpoint}\fi 
\begin{picture}(170,80)(0,0)

\put(-5,10){\vector(1,0){170.00}}

\put(0,5){\vector(0,1){74.00}}

\put(0.2,10){\line(0,1){22.20}}
\put(-0.2,10){\line(0,1){22.20}}
\put(0.4,10){\line(0,1){22.2}}
\put(-0.4,10){\line(0,1){22.20}}

\multiput(0,32)(4,0){41}{\line(1,0){2}}
\multiput(0,45)(4,0){41}{\line(1,0){2}}
\multiput(0,53)(4,0){41}{\line(1,0){2}}
\multiput(0,66)(4,0){41}{\line(1,0){2}}

\put(0,45){\circle*{1.0}}

\put(0.4,53){\line(0,1){13.20}}
\put(0.2,53){\line(0,1){13.20}}
\put(-0.2,53){\line(0,1){13.20}}
\put(-0.4,53){\line(0,1){13.20}}

\put(-9.5,11){$\scriptstyle\l_1^-=0$}
\put(-5.0,31){$\scriptstyle\l_1^+$}
\put(-5.0,65){$\scriptstyle\l_2^+$}
\put(-5.0,53){$\scriptstyle\l_2^-$}
\put(-4.0,7.0){$\scriptstyle 0=z_1^-$}
\put(23.0,6.8){$\scriptstyle z_1^+$}
\put(70.0,6.5){$\scriptstyle 2\pi-z_2^+$}
\put(98.0,6.5){$\scriptstyle 2\pi-z_1^+$}
\put(147.0,6.5){$\scriptstyle 2\pi+z_1^+$}
\put(12.0,11.5){$\scriptstyle \s_{1,0}(\O)$}
\put(102.0,11.5){$\scriptstyle \s_{1,1}(\O)$}
\put(139.0,11.5){$\scriptstyle \s_{1,2}(\O)$}
\put(40.5,11.5){$\scriptstyle \s_{2,0}(\Omega)$}
\put(75.5,11.5){$\scriptstyle \s_{2,1}(\Omega)$}
\put(17.0,59.0){$\l=1-\cos z$}
\put(-9.5,45){$\scriptstyle\s_3(\D)$}
\put(-9.5,59){$\scriptstyle\s_2(\D)$}
\put(-3.5,70){$\scriptstyle 2$}
\put(-9.5,20){$\scriptstyle\s_1(\D)$}
\put(-1,10){\line(1,0){2.00}}
\put(-1,70){\line(1,0){2.00}}

\multiput(-1,10)(4,0){41}{\line(1,0){2}}
\multiput(-1,70)(4,0){41}{\line(1,0){2}}
\bezier{600}(0,10)(15,9)(31.4,40)
\bezier{600}(94.2,40)(62.8,100)(31.4,40)
\bezier{600}(94.2,40)(125.6,-20)(157,40)

\put(34.2,10.0){\circle*{1.0}}
\put(28.0,7.0){$\scriptstyle\s_{3,0}(\O)$}
\put(86.0,7.){$\scriptstyle\s_{3,1}(\O)$}

\multiput(26.9,10.5)(0,2){11}{\line(0,1){1}}
\multiput(34.2,10.0)(0,2){18}{\line(0,1){1}}
\multiput(51.2,10.7)(0,2){28}{\line(0,1){1}}
\multiput(39.2,10.1)(0,2){22}{\line(0,1){1}}

\put(62.8,10.0){\circle*{1.0}}
\put(125.6,10.0){\circle*{1.5}}
\put(62.0,7.5){$\scriptstyle \pi$}
\put(124.0,6.5){$\scriptstyle 2\pi$}

\multiput(74.4,10.7)(0,2){28}{\line(0,1){1}}
\multiput(86.4,10.1)(0,2){22}{\line(0,1){1}}
\multiput(98.7,10.5)(0,2){11}{\line(0,1){1}}
\multiput(91.4,10.0)(0,2){18}{\line(0,1){1}}
\multiput(152.4,10.5)(0,2){11}{\line(0,1){1}}
\put(74.4,10.3){\line(1,0){12}}
\put(74.4,10.2){\line(1,0){12}}
\put(74.4,10.1){\line(1,0){12}}
\put(74.4,9.9){\line(1,0){12}}
\put(74.4,9.8){\line(1,0){12}}
\put(74.4,9.7){\line(1,0){12}}
\put(91.4,10.0){\circle*{1.0}}

\put(98.7,10.3){\line(1,0){53.8}}
\put(98.7,10.2){\line(1,0){53.8}}
\put(98.7,10.1){\line(1,0){53.8}}
\put(98.7,9.9){\line(1,0){53.8}}
\put(98.7,9.8){\line(1,0){53.8}}
\put(98.7,9.7){\line(1,0){53.8}}
\put(39.2,10.3){\line(1,0){12}}
\put(39.2,10.2){\line(1,0){12}}
\put(39.2,10.1){\line(1,0){12}}
\put(39.2,9.9){\line(1,0){12}}
\put(39.2,9.8){\line(1,0){12}}
\put(39.2,9.7){\line(1,0){12}}
\put(39.0,7.0){$\scriptstyle z_2^-$}
\put(50.0,7.0){$\scriptstyle z_2^+$}

\put(0,10.3){\line(1,0){26.9}}
\put(0,10.2){\line(1,0){26.9}}
\put(0,10.1){\line(1,0){26.9}}
\put(0,9.9){\line(1,0){26.9}}
\put(0,9.8){\line(1,0){26.9}}
\put(0,9.7){\line(1,0){26.9}}

\put(162.0,6.0){$z$}

\put(-4,76.0){$\l$}

\end{picture}
\caption{\footnotesize Relation between the spectra of $\D$ and $\O=\sqrt{\D_M}$.}
\label{fRel}
\end{figure}

\no \textbf{Remark.} 1) The relation between the spectra of $\D$ and $\O=\sqrt{\D_M}$ is shown in Fig.\ref{fRel}.

2) For each $j\in\N$ the spectral bands $\s_{1,2j-1}(\D_M)$ and $\s_{1,2j}(\D_M)$ defined in \er{sDM1} have the form
$$
\s_{1,2j-1}(\D_M)=[E_{1,2j-1}^-,E_{1,2j-1}^+]=[(2\pi j-z_1^+)^2,(2\pi j)^2],
$$
$$
\s_{1,2j}(\D_M)=[E_{1,2j}^-,E_{1,2j}^+]=[(2\pi j)^2,(2\pi j+z_1^+)^2].
$$
Thus, there is no gap between these spectral bands and the point $(2\pi j)^2$
is the point of tangency of these bands
(see Fig.\ref{fRel}).

 \subsection{Effective masses for metric Laplacians}

Let $\l_n(\vt)$, $n\in\N_{\n-r}$, $\vt\in\T^d$, be a band function
of the discrete Laplacian $\D$ and let $\l_n(\vt)$ have a minimum
(maximum) at some point $\vt_0$. Assume that $\l_n(\vt_0)$ is a
simple eigenvalue of $\D(\vt_0)$. We denote by $\m_n(\o)$ the effective form for the discrete Laplacian $\D$, defined in \er{lam}.
Due to \er{egv1}, the band
functions $E_{n,j}(\vt)$, $\vt\in\T^d$, $j=0,1,\ldots$\,, of the
metric Laplacian $\D_M$ have an extremum at the same point
$\vt_0$ and $E_{n,j}(\vt_0)$ are simple eigenvalues of $\D_M(\vt_0)$.
The analytic perturbation theory (see, e.g., \cite{K80}) gives that
the eigenvalues $E_{n,j}(\vt)$ have asymptotics:
\[\label{lam'}
E_{n,j}(\vt)=E_{n,j}(\vt_0)+\ve^2\m_{n,j}(\o)+O(\ve^3),\qqq
\m_{n,j}(\o)=\textstyle{1\/2}\,\ddot E_{n,j}(\vt_0+\ve\o)\big|_{\ve=0},
\]
as $\vt=\vt_0+\ve\o$, $\o\in \S^{d-1}$, $\ve=|\vt-\vt_0|\rightarrow0$, where
$\dot u=\pa u/\pa \ve$ and $\S^{d}$ is the $d$-dimensional sphere.

Now we formulate the relation between the effective forms for
discrete and metric Laplacians.

\begin{theorem}\lb{TEM}
i) Let a band function $E_{n,j}(\vt)$, $\vt\in\T^d$, have a minimum
(maximum) at some point $\vt_0$ and  some  ${n=1,2,\ldots,\n-r}$;
$j=0,1,\ldots$\,. Assume that $E_{n,j}(\vt_0)\neq(\pi k)^2$, $k\in\Z$,
is a simple eigenvalue of $\D_M(\vt_0)$. Then the effective form
$\mu_{n,j}(\o), \o\in \S^{d-1}$, for the metric Laplacian $\D_M$ at
the point $\vt_0$, defined in \er{lam'},  satisfies
\[
\lb{rel} \mu_{n,j}(\cdot)={2\sqrt{E_{n,j}(\vt_0)}\/\sin
\sqrt{E_{n,j}(\vt_0)}}\,\mu_n(\cdot),
\]
where $\m_n(\cdot)$ is the effective form for the discrete Laplacian at
the point $\vt_0$, defined in \er{lam}.

ii) The effective form $\mu_{1,0}(\o)$ at the beginning of the spectrum
 of the metric Laplacian $\D_M$  is non-degenerate and satisfies
\[\lb{mef}
0<{2\L_0\/\vk\,\n\,d}\leq\mu_{1,0}(\o)=2\mu(\o)\leq
{2\/\vk}\sum\limits_{\be\in\cB_*}\lan\t
(\be),\o\ran^2\leq{2\L_1\/\vk}\,, \qqq \vk=\sum\limits_{v\in V_*}\vk_v,
\]
for any $\o\in \S^{d-1}$, where $\mu(\o)$ is the effective form at
the beginning of the spectrum for the discrete Laplacian defined in
\er{lam}, $\n$ is the number of the
fundamental graph vertices, $\t (\be)$ is the index of the edge
$\be\in\cE_\ast$; $\L_0$, $\L_1$ are defined before Theorem
\ref{Tem}.

\end{theorem}

\no {\bf Proof.} i) Due to \er{egv1}, we have
\[
\lb{EEE}
\cos z_{n,j}(\vt)=1-\l_n(\vt),\qq \textrm{where} \qq z_{n,j}(\vt)=\sqrt{E_{n,j}(\vt)}\;.
\]
Differentiating \er{EEE} with respect to $\ve$ twice and using that $\dot E_{n,j}(\vt_0+\ve\o)\big|_{\ve=0}=0$, we obtain
\[
\begin{aligned}
{\sin z_{n,j}(\vt_0+\ve\o)\/2\,z_{n,j}(\vt_0+\ve\o)}\,\dot E_{n,j}(\vt_0+\ve\o)=
\dot\l_n(\vt_0+\ve\o),
\end{aligned}
\]
\[
\lb{d2E}
{\sin z_{n,j}(\vt_0)\/2\, z_{n,j}(\vt_0)}\ddot E_{n,j}(\vt_0+\ve\o)\big|_{\ve=0}=
\ddot\l_n(\vt_0+\ve\o)\big|_{\ve=0}.
\]
Using the definitions of the effective forms $\mu_n(\o)$ and $\m_{n,j}(\o)$ in \er{lam} and \er{lam'}, we have the identity \er{rel}.

ii) Due to \er{egv1}, we have
\[
\lb{EEE1}
\cos z_{1,0}(\vt)=1-\l_1(\vt),\qq \textrm{where} \qq z_{1,0}(\vt)=\sqrt{E_{1,0}(\vt)}\;.
\]
The eigenvalues $\l_1(\vt)$ and $E_{1,0}(\vt)$ have asymptotics as $\ve=|\vt|\rightarrow0$:
\[\label{lam1}
\l_1(\vt)=\ve^2\m(\o)+O(\ve^3), \qqq E_{1,0}(\vt)=\ve^2\m_{1,0}(\o)+O(\ve^3).
\]
Substituting \er{lam1} into \er{EEE1} and expanding the left side of \er{EEE1} into a Taylor series with respect to $\ve$, we obtain
\[
\lb{EEE2}
1-\ve^2\,{\m_{1,0}(\o)\/2}+O(\ve^3)=1-\ve^2\m(\o)+O(\ve^3),
\]
which yields $\m_{1,0}(\o)=2\m(\o)$. This and \er{tseN} give \er{mef}. \qq \BBox

\section{\lb{Sec4} Effective masses for combinatorial Laplacians}
\setcounter{equation}{0}

The combinatorial Laplacian $\wh\D$ on $f\in\ell^2(V)$ is defined by

\[
\lb{DOL} \big(\wh\D f\big)(v)=
\sum\limits_{(v,\,u)_e\in\cE}\big(f(v)-f(u)\big), \qqq
 v\in V.
\]

The discrete Laplacian $\wh \D$ on $\ell^2(V)$ has the standard decomposition into a
constant fiber direct integral  by
\[
\lb{razN}
\ell^2(V)={1\/(2\pi)^d}\int^\oplus_{\T^d}\ell^2(V_\ast)\,d\vt ,\qqq
U\wh \D U^{-1}={1\/(2\pi)^d}\int^\oplus_{\T^d}\wh \D(\vt)d\vt,
\]
$\T^d=\R^d/(2\pi\Z)^d$,
for some unitary operator $U$. Here $\ell^2(V_\ast)=\C^\nu$ is the fiber space
and $\wh \D(\vt)=\{\wh \D_{uv}(\vt)\}_{u,v\in V_\ast}$ is the Floquet $\nu\ts\nu$  matrix, having the form
\[
\lb{l2.15'}
\wh\D_{uv}(\vt)=\vk_u\d_{uv}-\ca \sum\limits_{{\bf
e}=(u,\,v)\in{\cA}_\ast}e^{\,i\lan\t
(\be),\,\vt\ran }, \qq &  {\rm if}\  \ (u,v)\in \cA_\ast \\
\qqq 0, &  {\rm if}\  \ (u,v)\notin \cA_\ast \ac,
\]
see \cite{KS14}, where $\lan\cdot\,,\cdot\ran$ denotes the standard inner
product in $\R^d$.

Let $\wh\l(\vt)$, $\vt\in\T^d$, be a band function of the Laplacian
$\wh\D$ and let $\wh\l(\vt)$ have a minimum (maximum) at some point
$\vt_0$. Assume that $\wh\l(\vt_0)$ is a simple eigenvalue of
$\wh\D(\vt_0)$. Then the eigenvalue $\wh\l(\vt)$ has
the asymptotics as $\vt=\vt_0+\ve\o$, $\o\in\S^{d-1}$, $\ve=|\vt-\vt_0|\to0$:
\[\label{lamC}
\begin{aligned}
\wh\l(\vt)=\wh\l(\vt_0)+\ve^2\wh\m(\o)+O(\ve^3), \qqq
\textstyle\wh\m(\o)={1\/2}\,\ddot{\wh\l}(\vt_0+\ve\o)\big|_{\ve=0},\qq
\end{aligned}
\]
where $\dot u=\pa u/\pa \ve$ and $\S^{d}$ is the $d$-dimensional
sphere.

\begin{theorem}\lb{T2N}
Let a band function $\wh\l(\vt)$, $\vt\in\T^d$, have a minimum
(maximum) at some point $\vt_0$. Assume that $\wh\l(\vt_0)$ is a
simple eigenvalue of $\wh\D(\vt_0)$. Then the effective form
$\wh\mu(\o)$ from \er{lamC} satisfies
\[
\lb{em11''} \big|\wh\mu(\o)\big|\leq
{T_1^2\/\wh\r}+T_2,\qqq\text{where}\qq T_s={1\/s}\,\max_{u\in
V_\ast}\sum_{\be=(u,v)\in\cA_\ast} \|\t (\be)\|^s,\qq s=1,2,
\]
$\wh\r=\wh\r(\vt_0)$ is the distance between $\wh\l(\vt_0)$  and the set
$\s\big(\wh\D(\vt_0)\big)\setminus\big\{\wh\l(\vt_0)\big\}$, $\t
(\be)$ is the index of the edge $\be\in\cA_\ast$.
\end{theorem}

Recall that we orient the undirected
edges of the fundamental graph $\G_*$ from the set $\cE_*$ in some arbitrary way and denote by $\cB_*$ the set of all bridges from $\cE_*$.

\

\begin{proposition}\lb{TP1C}
Each Floquet matrix $\wh \D(\vt)$, $\vt\in\T^d$, has the following factorization
\[\lb{DvtC}
\wh \D(\vt)=\wh \na^\ast(\vt)\wh \na(\vt),
\]
where the $\n_\ast\ts\n$ matrix $\wh \na(\vt)=\na(\vt)\vk^{1\/2}$, $\vk=\diag(\vk_v)_{v\in V_*}$, $\na(\vt)$ is defined by \er{Nvt}.
\end{proposition}

We describe the effective masses for the case of the beginning  the
spectrum.

\begin{theorem}
\label{TemC}
Let $\wh\mu(\o)$ be the effective form at the beginning of the spectrum defined in \er{lamC}. Then each $\wh\mu(\o)$, $\o\in\S^{d-1}$, satisfies
\[
\lb{em11.1C}
\textstyle \wh\mu(\o)=\|P_*\wh h(\o)\|_*^2,
\]
\[
\lb{tseC}
0<{\L_0\/\n^2\, d}\leq\wh\mu(\o)\leq\|\wh h(\o)\|^2_*=
{1\/\n}\sum\limits_{\be\in\cB_*}\lan\t
(\be),\o\ran^2\leq{\L_1\/\n}\,,
\]
where $P_*$ is the projection onto the kernel of $\wh\na^*(0)$ and
the vector $\wh h(\o)$ is given by
\[\lb{hhhC}
\wh h(\o)=
{i\/\sqrt{\n}}\big(\lan\t(\be),\o\ran\big)_{\be\in\cE_\ast}\in\C^{\n_\ast},
\]
$\n$ is the number of the fundamental graph vertices, $\n_\ast$ is the number of the fundamental graph edges from $\cE_\ast$,
$\|\cdot\|_\ast$ denotes the norm in $\C^{\n_\ast}$, $\t (\be)$ is
the index of the edge $\be\in\cE_\ast$; $\L_0$, $\L_1$ are defined before Theorem \ref{Tem}.
\end{theorem}

Let a  matrix $\wh M\ge 0$ be defined by the effective form
$\wh\mu(\o)=\lan\wh M\o,\o\ran$. Then  the effective mass tensor has the form
$
\wh m=\wh M^{-1}.
$
Let  $\wh m_1\ge \wh m_2\ge\ldots\ge \wh m_d$ be eigenvalues of $\wh
m$.

\begin{corollary}\label{TsemC}
Let  $\be_0\in\cB_\ast$.   Then at the beginning  of
the spectrum  the following estimates hold true
\[
\lb{tem1x} \wh C_0\leq\Tr \wh M\leq \wh C_1,
\]
\[
\lb{tem2x} {\n^2\,d\/\L_0}\ge \wh m_1\ge\ldots\ge \wh m_d\ge {1\/\wh
C_1}\,,
\]
where
\[
\lb{tem3x} \wh C_0={1\/\n^2}\,\|\t(\be_0)\|^2,\qqq \wh C_1={1\/\n}
\sum\limits_{\be\in\cB_\ast}\|\t (\be)\|^2.
\]
\end{corollary}

We omit the proof of these results since the proof for the
combinatorial  Laplacian repeats the proof for the normalized  one.

\section{\lb{Sec3x} Examples of effective masses}
\setcounter{equation}{0}

\subsection{Effective masses for the lattice graph.}
We consider the lattice graph $\dL^d=(V,\cE)$, where the vertex set  and the edge set are given by
\[
\lb{dLg}
V=\Z^d,\qquad  \cE=\big\{(m,m+a_1),\ldots,(m,m+a_d), \quad
\forall\,m\in\Z^d\big\},
\]
and $a_1,\ldots,a_d$ is the standard orthonormal basis, see
Fig.\ref{ff.0.1}\emph{a}. The
"minimal"\, fundamental graph $\dL_\ast^d$ of the lattice $\dL^d$
consists of one vertex $v$ and $d$ edge-loops
$(v,v)$, see Fig.\ref{ff.0.1}\emph{b}. All
bridges of $\dL_\ast^d$ are loops and their indices have
the form $a_1,\ldots,a_d$.

\setlength{\unitlength}{1.0mm}
\begin{figure}[h]
\centering

\unitlength 1.0mm 
\linethickness{0.4pt}
\ifx\plotpoint\undefined\newsavebox{\plotpoint}\fi 
\begin{picture}(140,60)(0,0)


\put(10,10){\line(1,0){40.00}}
\put(10,30){\line(1,0){40.00}}
\put(10,50){\line(1,0){40.00}}

\bezier{60}(17,33.5)(37,33.5)(57,33.5)
\bezier{60}(17,13.5)(37,13.5)(57,13.5)

\put(17,53.5){\line(1,0){40.00}}
\put(10,10){\line(0,1){40.00}}
\put(30,10){\line(0,1){40.00}}
\put(50,10){\line(0,1){40.00}}
\put(57,13.5){\line(0,1){40.00}}
\bezier{12}(10,10)(13.5,11.75)(17,13.5)

\bezier{60}(17,13.5)(17,33.5)(17,53.5)
\bezier{60}(37,13.5)(37,33.5)(37,53.5)
\bezier{12}(30,10)(33.5,11.75)(37,13.5)
\put(50,10){\line(2,1){7.00}}
\bezier{12}(10,30)(13.5,31.75)(17,33.5)
\bezier{12}(30,30)(33.5,31.75)(37,33.5)

\put(50,30){\line(2,1){7.00}}
\put(10,50){\line(2,1){7.00}}
\put(30,50){\line(2,1){7.00}}
\put(50,50){\line(2,1){7.00}}
\put(10,10){\circle{1.0}}
\put(30,10){\circle{1.0}}
\put(50,10){\circle*{1.0}}
\put(10,30){\circle{1.0}}
\put(30,30){\circle{1.0}}
\put(50,30){\circle{1.0}}
\put(10,50){\circle{1.0}}
\put(30,50){\circle{1.0}}
\put(50,50){\circle{1.0}}

\put(17,53.5){\circle{1.0}}
\put(37,53.5){\circle{1.0}}
\put(57,53.5){\circle{1.0}}
\put(17,33.5){\circle{1.0}}
\put(37,33.5){\circle{1.0}}
\put(57,33.5){\circle{1.0}}
\put(17,13.5){\circle{1.0}}
\put(37,13.5){\circle{1.0}}
\put(57,13.5){\circle{1.0}}

\put(50,10){\vector(0,1){20.00}}
\put(50,10){\vector(-1,0){20.00}}
\put(50,10){\vector(2,1){7.00}}

\put(49.0,7.0){$v$}
\put(54,10.0){$a_1$}
\put(35.0,7.5){$a_2$}
\put(46.0,26){$a_3$}

\put(-5,5){(\emph{a})}
\bezier{30}(87,13.5)(97,13.5)(107,13.5)
\bezier{30}(87,13.5)(87,23.5)(87,33.5)
\bezier{12}(87,13.5)(83.5,11.75)(80,10)

\put(100,10){\vector(2,1){7.00}}
\bezier{12}(100,30)(103.5,31.75)(107,33.5)
\bezier{12}(80,30)(83.5,31.75)(87,33.5)
\bezier{30}(80,10)(80,20)(80,30)
\bezier{30}(107,13.5)(107,23.5)(107,33.5)
\put(100,10){\vector(-1,0){20.00}}
\bezier{30}(87,33.5)(97,33.5)(107,33.5)

\put(100,10){\vector(0,1){20.00}}
\bezier{30}(80,30)(90,30)(100,30)
\put(80,10){\circle{1}}
\put(80,30){\circle{1}}
\put(100,30){\circle{1}}
\put(87,13.5){\circle{1}}
\put(107,13.5){\circle{1}}
\put(107,33.5){\circle{1}}
\put(87,33.5){\circle{1}}

\put(100,10){\circle*{1}}

\put(100.0,7.0){$v$}
\put(79.0,7.0){$v$}
\put(101.0,28.5){$v$}
\put(77.0,28.5){$v$}

\put(84.5,33.5){$v$}
\put(108.0,33.5){$v$}
\put(108.0,13.5){$v$}

\put(103.0,9.5){$\be_1$}
\put(89.0,7.0){$\be_2$}
\put(96.0,21.0){$\be_3$}

\put(65,5){(\emph{b})}
\end{picture}

\vspace{-0.5cm} \caption{\footnotesize \emph{a}) Lattice
$\dL^3$;\quad \emph{b})  the fundamental graph $\dL^3_\ast$.}
\label{ff.0.1}
\end{figure}

For the lattice graph we have
$$
\l(\vt)=1-{1\/d}\,(\cos\vt_1+\ldots+\cos\vt_d)=\ve^2\m^-(\o)+O(\ve^3), \qqq \textrm{as } \ \ve=|\vt|\rightarrow0,
$$
where the effective form $\m^-(\o)$ associated with the beginning of the spectrum is given by
$$
\m^-(\o)={1\/2d}\,(\o_1^2+\ldots+\o_d^2)={1\/2d}\,,\qqq \o=(\o_s)_{s\in\N_d}={\vt\/|\vt|}\in\S^{d-1}\,.
$$
Since $\n=1$, $\vk=2d$ and the indices of all bridges form an orthonormal basis in $\R^d$,
$\L_0=\L_1=1$ and the estimate \er{tseN} yields
\[
{1\/2d^2}\leq\mu^-(\o)\leq{1\/2d}\,.
\]
Thus, the upper estimate in \er{tseN} is sharp.
The lattice graph is bipartite. Then the effective form associated with the  upper end of the spectrum $\m^+(\o)=-\m^-(\o)$. We note that the effective mass for the Laplacian on the lattice graph is independent of direction.

\subsection{Effective masses for the graphene.} We consider the hexagonal lattice $\bG$ (see Fig.\ref{ff.0.3}\emph{a}). The lattice $\bG$ is invariant
under translations through the vectors $a_1$, $a_2$.
The fundamental vertex set $V_0=\{v_1,v_2\}$ and the fundamental graph $\bG_*$ are also shown in the figure.
The fundamental graph $\bG_\ast$ consists of two vertices $v_1,v_2$ and \emph{multiple} edges $\be_1=\be_2=\be_3=(v_1,v_2)$
with the indices
$\t(\be_1)=(0,0)$, $\t(\be_2)=(1,0)$, $\t(\be_3)=(0,1)$ (Fig.\ref{ff.0.3}\emph{b}).

The Floquet matrix $\D(\vt)$ defined by \er{l2.15N} for
the hexagonal lattice $\bG$ has the form
\[
\lb{ell}
\D(\vt)=\left(
\begin{array}{cc}
1 & -\D_{12}(\vt) \\[6pt]
-\bar\D_{12}(\vt) & 1
\end{array}\right),\qqq \textstyle \D_{12}(\vt)={1\/3}\,\big(1+e^{i\vt _1}+
e^{i\vt_2}\big),
\]
for all $\vt=(\vt_1,\vt_2)\in\T^2$.
 This yields that the eigenvalues of each matrix $\D(\vt)$
are given by
\[
\lb{evg}\textstyle \l_n(\vt)=1+(-1)^{n}|\D_{12}(\vt
)|\;,  \qqq n=1,2.
\]

The band function $\l_1(\vt)$ attains its minimum only at the point 0 and
attains its maximum only at the points $\pm\vt_0$,
\[\lb{vt0}
\textstyle\vt_0=(\vt^0_1,\vt^0_2)=
\big(\frac{2\pi}3\,,-\frac{2\pi}3\big).
\]
The band function $\l_2(\vt)$ attains its maximum only at the point 0 and
attains its minimum only at the points $\pm\vt_0$.

Expanding the band function $\l_1(\vt)$ in the Taylor series about
the point $0$ and using that the hexagonal lattice is bipartite, we obtain
\[
\textstyle\m_1^-(\o)={1\/9}\,\big(\o_1^2+\o_2^2-\o_1\o_2\big),\qqq
\m_2^+(\o)=-\m_1^-(\o),\qqq \forall\, \o=(\o_1,\o_2)\in\S^1,
\]
where $\m_n^-(\o)$ and $\m_n^+(\o)$, $n=1,2$, are the effective forms associated with the lower and upper ends of the spectral band $\s_n(\D)$, respectively.

Since the matrix of the effective form $\m_1^-(\o)$ has the eigenvalues ${1\/18}$ and ${1\/6}$\,, we have the following estimates
\[
\textstyle{1\/18}\leq\m_1^-(\o)\leq{1\/6}\,.
\]
On the other hand, the fundamental graph $\bG_*$ has only two bridges with indices $(1,0)$, $(0,1)$, $d=2$, $\vk=6$, $\n=2$. Then $\L_0=\L_1=1$ and the estimate \er{tseN} gives
\[
\textstyle{1\/24}\leq\mu_1^-(\o)\leq{1\/6}\,.
\]
Each of the matrices $\D(\vt_0)$ and $\D(-\vt_0)$, where $\vt_0$ is defined in \er{vt0}, has the eigenvalue 1 of multiplicity 2. Therefore, there is no gap between the spectral bands.

\subsection{Effective masses for the stanene lattice}
Stanene is a topological insulator, theoretically predicted by
Shoucheng Zhang's group at Stanford, which may display
dissipationless currents at its edges near room temperature
\cite{Z13}. It is composed of tin atoms arranged in a single layer,
in a manner similar to graphene. Stanene got its name by combining
stannum (the Latin name for tin) with the suffix -ene used by
graphene. The addition of fluorine atoms to the tin lattice could
extend the critical temperature up to $100^\circ$ C. This would make
it practical for use in integrated circuits to make smaller, faster
and more energy efficient computers. Stanene has a band gap, it is a
semiconducting material. That makes it useful as material for use in
a transistor, which must have a component that turns on and off. For
more details see \cite{Z13} and references therein.

\setlength{\unitlength}{1.0mm}
\begin{figure}[h]
\centering

\unitlength 1mm 
\linethickness{0.4pt}
\ifx\plotpoint\undefined\newsavebox{\plotpoint}\fi 
\begin{picture}(100,45)(0,0)

\put(3,5){(\emph{a})}

\put(14,10){\circle{1}}
\put(14,10){\line(0,-1){4.00}}
\put(14,6){\circle{1}}
\put(14,22){\line(0,-1){4.00}}
\put(14,18){\circle{1}}
\put(14,34){\line(0,-1){4.00}}
\put(14,30){\circle{1}}

\put(28,10){\circle{1}}
\put(28,10){\line(0,1){4.00}}
\put(28,14){\circle{1}}
\put(28,22){\line(0,1){4.00}}
\put(28,26){\circle*{1}}
\put(28,34){\line(0,1){4.00}}
\put(28,38){\circle{1}}

\put(34,10){\circle{1}}
\put(34,10){\line(0,-1){4.00}}
\put(34,6){\circle{1}}
\put(34,22){\line(0,-1){4.00}}
\put(34,18){\circle{1}}
\put(34,34){\line(0,-1){4.00}}
\put(34,30){\circle{1}}

\put(48,10){\circle{1}}
\put(48,10){\line(0,1){4.00}}
\put(48,14){\circle{1}}
\put(48,22){\line(0,1){4.00}}
\put(48,26){\circle{1}}
\put(48,34){\line(0,1){4.00}}
\put(48,38){\circle{1}}

\put(18,16){\circle{1}}
\put(18,16){\line(0,1){4.00}}
\put(18,20){\circle{1}}
\put(18,28){\line(0,1){4.00}}
\put(18,32){\circle{1}}
\put(18,40){\line(0,1){4.00}}
\put(18,44){\circle{1}}

\put(24,16){\circle*{1}}
\put(24,16){\line(0,-1){4.00}}
\put(24,12){\circle*{1}}
\put(24,28){\line(0,-1){4.00}}
\put(24,24){\circle{1}}
\put(24,40){\line(0,-1){4.00}}
\put(24,36){\circle{1}}

\put(38,16){\circle{1}}
\put(38,16){\line(0,1){4.00}}
\put(38,20){\circle{1}}
\put(38,28){\line(0,1){4.00}}
\put(38,32){\circle{1}}
\put(38,40){\line(0,1){4.00}}
\put(38,44){\circle{1}}

\put(44,16){\circle{1}}
\put(44,16){\line(0,-1){4.00}}
\put(44,12){\circle{1}}
\put(44,28){\line(0,-1){4.00}}
\put(44,24){\circle{1}}
\put(44,40){\line(0,-1){4.00}}
\put(44,36){\circle{1}}

\put(14,22){\circle{1}}
\put(28,22){\circle*{1}}

\put(34,22){\circle{1}}
\put(48,22){\circle{1}}

\put(18,28){\circle{1}}
\put(24,28){\circle{1}}
\put(38,28){\circle{1}}
\put(44,28){\circle{1}}

\put(14,34){\circle{1}}
\put(28,34){\circle{1}}
\put(34,34){\circle{1}}
\put(48,34){\circle{1}}

\put(18,40){\circle{1}}
\put(24,40){\circle{1}}
\put(38,40){\circle{1}}
\put(44,40){\circle{1}}

\put(28,10){\line(1,0){6.00}}
\put(18,16){\line(1,0){6.00}}
\put(38,16){\line(1,0){6.00}}

\put(28,22){\line(1,0){6.00}}
\put(18,28){\line(1,0){6.00}}
\put(38,28){\line(1,0){6.00}}

\put(28,34){\line(1,0){6.00}}
\put(18,40){\line(1,0){6.00}}
\put(38,40){\line(1,0){6.00}}

\put(14,10){\line(2,3){4.00}}
\put(34,10){\line(2,3){4.00}}
\put(24,16){\line(2,3){4.00}}
\put(44,16){\line(2,3){4.00}}

\put(14,22){\line(2,3){4.00}}
\put(34,22){\line(2,3){4.00}}
\put(24,28){\line(2,3){4.00}}
\put(44,28){\line(2,3){4.00}}

\put(14,34){\line(2,3){4.00}}
\put(34,34){\line(2,3){4.00}}

\put(28,10){\line(-2,3){4.00}}
\put(48,10){\line(-2,3){4.00}}
\put(38,16){\line(-2,3){4.00}}
\put(18,16){\line(-2,3){4.00}}

\put(28,22){\line(-2,3){4.00}}
\put(48,22){\line(-2,3){4.00}}
\put(38,28){\line(-2,3){4.00}}
\put(18,28){\line(-2,3){4.00}}

\put(28,34){\line(-2,3){4.00}}
\put(48,34){\line(-2,3){4.00}}

\put(27.8,19.5){$\scriptstyle v_1$}
\put(28.8,26.5){$\scriptstyle v_3$}
\put(22.0,17.5){$\scriptstyle v_2$}
\put(21.0,10.5){$\scriptstyle v_4$}

\put(75,10){\circle*{1}}
\put(75,10){\line(0,-1){6.00}}
\put(75,4){\circle*{1}}

\put(83,21){\circle*{1}}
\put(83,21){\line(0,1){6.00}}
\put(83,27){\circle*{1}}

\put(75,30){\circle{1}}
\put(95,40){\circle{1}}
\put(95,20){\circle{1}}

\put(75,10){\vector(0,1){20.0}}
\put(75,10){\vector(2,1){20.0}}

\multiput(95,20)(0,7){3}{\line(0,1){4}}
\put(75,30){\line(2,1){4.0}}
\put(82,33.5){\line(2,1){4.0}}
\put(89,37){\line(2,1){4.0}}

\qbezier(83,21)(89,20.5)(95,20)
\qbezier(83,21)(79,15.5)(75,10)
\qbezier(83,21)(79,25.5)(75,30)

\put(70.5,9.0){$v_2$}
\put(70.5,4.0){$v_4$}
\put(83,18){$v_1$}
\put(83,28){$v_3$}
\put(96,19.0){$v_2$}
\put(71.0,31.0){$v_2$}
\put(93.5,42.0){$v_2$}
\put(85,12.5){$a_1$}
\put(70.5,20.0){$a_2$}

\put(77,17.2){$\scriptstyle\be_1$}
\put(89,21.2){$\scriptstyle\be_2$}
\put(78,27.0){$\scriptstyle\be_3$}
\put(83.5,23.5){$\scriptstyle\be_4$}
\put(75.5,6.2){$\scriptstyle\be_5$}

\put(62,5){(\emph{b})}
\end{picture}

\vspace{-0.5cm} \caption{\footnotesize  \emph{a}) Stanene lattice
$\bS$; \quad \emph{b}) the fundamental graph $\bS_\ast$.} \label{ff.FCC}
\end{figure}

The stanene lattice
$\bS$ is obtained from the hexagonal lattice $\bG$ by adding a pendant edge at each vertex of $\bG$ (Fig.\ref{ff.FCC}\emph{a}).
The fundamental graph $\bS_\ast$ of $\bS$ consists of 4 vertices
$v_1,v_2,v_3,v_4$ and 5 edges (Fig.\ref{ff.FCC}\emph{b})
$$
\be_1=\be_2=\be_3=(v_1,v_2),\qqq \be_4=(v_1,v_3),\qqq \be_5=(v_2,v_4).
$$
The indices of the fundamental graph edges with respect to the fundamental vertex set $V_0=\{v_1,v_2,v_3,v_4\}$, see Fig.\ref{ff.FCC}\emph{a}, are given by
$$
\t(\be_2)=(1,0),\qqq \t(\be_3)=(0,1),\qqq
\t(\be_1)=\t(\be_4)=\t(\be_5)=(0,0).
$$
For each $\vt\in\T^2$ the Floquet matrix $\D(\vt)$ defined by \er{l2.15N}
has the form
$$
\D(\vt)=\1_4-
{1\/2}\begin{pmatrix}
0&b(\vt)&1&0\\
\overline b(\vt)&0&0&1 \\
1&0&0&0 \\
0&1&0&0 \\
\end{pmatrix},\qqq b(\vt)={1\/2}\big(1+e^{i\vt_1}+e^{i\vt_2}\big).
$$
We have the characteristic equation for the matrix $\D(\vt)$
$$
(\l-1)^4-(\l-1)^2\left(\frac{1}{2}+\frac{|b(\vt)|^2}{4}\right)
+\frac{1}{16}=0.
$$
The eigenvalues of each matrix $\D(\vt)$ are given by
$$
\textstyle \l_{1,2,3,4}(\vt) =1 \pm\frac{|b(\vt)|}4\pm\frac{\sqrt{|b(\vt)|^2+4}}{4}\,.
$$

The band function
$$
\textstyle\l_{1}(\vt) =1-\frac{|b(\vt)|}4-\frac{\sqrt{|b(\vt)|^2+4}}{4}
$$
attains its minimum only at the point 0 and
attains its maximum only at the points $\pm\vt_0$,
\[\lb{Svt0}
\textstyle\vt_0=(\vt^0_1,\vt^0_2)=
\big(\frac{2\pi}3\,,-\frac{2\pi}3\big).
\]
The band function
$$
\textstyle\l_{2}(\vt) =1+\frac{|b(\vt)|}4-\frac{\sqrt{|b(\vt)|^2+4}}{4}
$$
attains its maximum only at the point 0 and
attains its minimum only at the points $\pm\vt_0$.

Expanding the band functions $\l_1(\vt)$ and $\l_2(\vt)$ in the Taylor series about the point $0$ and using that the stanene lattice is bipartite, we obtain
\[
\textstyle\m_1^-(\o)={1\/15}\,\big(\o_1^2+\o_2^2-\o_1\o_2\big),\qqq
\m_4^+(\o)=-\m_1^-(\o)=-{1\/15}\,\big(\o_1^2+\o_2^2-\o_1\o_2\big),
\]
\[
\textstyle\m_2^+(\o)=-{1\/60}\,\big(\o_1^2+\o_2^2-\o_1\o_2\big),\qqq
\m_3^-(\o)=-\m_2^+(\o)={1\/60}\,\big(\o_1^2+\o_2^2-\o_1\o_2\big).
\]
Since the matrix of the effective form $\m_1^-(\o)$ has the eigenvalues ${1\/30}$ and ${1\/10}$ and the matrix of $\m_2^+(\o)$ has the eigenvalues $-{1\/40}$ and $-{1\/120}$\,, we have the following estimates
\[
\textstyle{1\/30}\leq\m_1^-(\o)\leq{1\/10}\,,\qqq
-{1\/40}\leq\m_2^+(\o)\leq-{1\/120}\,.
\]
On the other hand, the fundamental graph $\bS_*$ has only two bridges with indices $(1,0)$, $(0,1)$ forming an orthonormal basis in $\R^2$, $d=2$, $\vk=\sum\limits_{j=1}^4\vk_{v_j}=10$, $\n=4$. Then $\L_0=\L_1=1$ and the estimate \er{tseN} gives
\[
\textstyle{1\/80}={1\/10\cdot4\cdot2}\leq\m_1^-(\o)\leq{1\/\vk}={1\/10}\,.
\]

Each of the matrices $\D(\vt_0)$ and $\D(-\vt_0)$, where $\vt_0$ is defined
in \er{Svt0}, has the eigenvalues ${1\/2}$ and ${3\/2}$ of multiplicity 2. Therefore, there is no gap between the spectral bands $\s_1(\D)$, $\s_2(\D)$ and between the spectral bands $\s_3(\D)$, $\s_4(\D)$.
Note that the spectrum of the Laplacian on $\bS$ has the form
$$
\textstyle\s(\D)=\s_{ac}(\D)=\big[0;{3\/4}\,\big]\cup\big[{5\/4}\,;2\big].
$$

\medskip

\setlength{\itemsep}{-\parskip} \footnotesize \no {\bf
Acknowledgments.} { Evgeny Korotyaev's  study was partly supported
by the RFFI grant  No 11-01-00458 and by the  project  SPbGU No
11.38.215.2014.}

\end{document}